\documentclass[11pt]{amsart}
\usepackage[paperwidth=6.7in,paperheight=9.6in,inner=0.85in,outer=0.85in,top=0.9in,bottom=0.9in]{geometry}
\usepackage[T1]{fontenc}
\usepackage{lmodern,microtype}
\usepackage{amsmath,amssymb,amsthm,mathtools,mathrsfs}
\usepackage[hidelinks]{hyperref}
\hypersetup{pdftitle={Generic Strong Regularity and the Branson\textendash Gover Conjecture},
pdfsubject={Generic surjectivity of conformal harmonics and strong regularity of Branson-Gover operators},
pdfkeywords={Branson-Gover conjecture, generic surjectivity, conformal harmonics, Branson-Gover operators, conformal Hodge theory, metric variations},
pdfauthor={Wei-Chun Chen, Wei-Ting Kao, Chun-Kai Lien}}
\newtheorem{theorem}{Theorem}[section]
\newtheorem{proposition}[theorem]{Proposition}
\newtheorem{lemma}[theorem]{Lemma}
\newtheorem{corollary}[theorem]{Corollary}
\newcommand{\Om}{\Omega}
\newcommand{\HH}{\mathcal H}
\newcommand{\Met}{\operatorname{Met}}
\newcommand{\im}{\operatorname{im}}
\newcommand{\Id}{\operatorname{Id}}
\newcommand{\Sym}{\operatorname{Sym}}
\allowdisplaybreaks[1]
\title{Generic Strong Regularity\\and the Branson--Gover Conjecture}
\author{Wei-Chun Chen}
\address{Department of Mathematics, National Taiwan University, Taipei 10617, Taiwan}
\email{b14201028@ntu.edu.tw}

\author{Wei-Ting Kao}
\address{Institute of Mathematics, Academia Sinica, Taipei 106319, Taiwan}
\email{wtkao@as.edu.tw}

\author{Chun-Kai Lien}
\address{Department of Mathematics, National Taiwan University, Taipei 10617, Taiwan}
\email{d13221002@ntu.edu.tw}
\date{}
\subjclass[2020]{53C18, 58J50, 58A10}
\keywords{Conformal geometry, Branson--Gover operators,
conformal harmonics, generic metrics, metric variations}
\begin{document}

\begin{abstract}
On each closed smooth manifold of even dimension \(2N\ge 4\), we prove that the Branson--Gover detour operators are strongly regular in every degree \(k-1\), \(1\le k<N\), on a \(C^\infty\)-open dense set of Riemannian metrics. Consequently, every de Rham cohomology class of positive degree below the middle dimension has a unique conformal harmonic representative for these metrics. On a connected manifold, the critical scalar GJMS operator has only constants in its kernel on an open dense set of conformal structures modulo diffeomorphisms isotopic to the identity, as conjectured by Canzani, Gover, Jakobson, and Ponge. The proof computes the principal symbol of a reduced metric Hessian that accounts for the motion of the kernel through a Green operator.
\end{abstract}

\maketitle
\section{Introduction}\label{sec:introduction}

Branson and Gover introduced conformal harmonic forms and proposed
that, generically, they represent every de Rham class of positive degree
below the middle dimension \cite[discussion after Corollary 2.13]{BG}.
Case and Gover state this as the surjectivity of the natural map to
de Rham cohomology \cite[Conjecture 3.3]{CG}. Gover \cite[p.~122]{Gover2005}
and Branson \cite[p.~43]{Branson2005} proposed that detour cohomology
generically agrees with de Rham cohomology. We prove strong regularity
on an open dense set of metrics on each closed even-dimensional
manifold. Strong regularity implies the conjectured surjectivity.

Let \(M\) be a closed smooth manifold of dimension \(2N\ge4\).
{Write \(\Om^r=\Om^r(M)\) for the smooth \(r\)-forms.}
We consider real forms and denote by \(\delta_g\) the formal adjoint
of \(d\) with respect to the Riemannian density. The Hodge Laplacian
\(\Delta_g=d\delta_g+\delta_gd\) is nonnegative. We write \((u,v)_g\)
for the integrated \(L^2_g\) pairing and \(\langle u,v\rangle_{g,x}\)
for the pointwise pairing. For the symbol calculations we complexify
the bundles.

\begin{samepage}
Fix \(1\le k<N\). Let \(Q_{k,g}\) be the canonical Branson--Gover
operator on closed \(k\)-forms, defined and normalized in
Section~\ref{subsec:canonical-definition} by \eqref{eq:ambient-definition},
with weighted targets identified with ordinary forms using \(g\). Set
\(L_g=\delta_gQ_{k,g}d\).
\end{samepage}
The ambient construction in Section~\ref{subsec:canonical-definition}
gives locality, naturality and smooth dependence on finitely many metric
jets. The principal parts of \(Q_{k,g}\) on closed forms and of \(L_g\)
are \((d\delta_g)^{N-k}\) and \((\delta_gd)^{N-k+1}\), respectively.
The integrated pairing of \(Q_{k,g}\) on closed forms is symmetric
\cite[Theorem 2.8(i)]{BG}. For a constant \(c>0\), the correction in the
conformal transformation law vanishes on closed forms
\cite[Theorem 2.8(iv)]{BG}. The output weight then gives
\(Q_{k,c^2g}=c^{2k-2N}Q_{k,g}\).

We say that \(g\) is strongly regular in degree \(k-1\) if
\(\ker L_g=\ker d\).

\begin{theorem}[Generic strong regularity]\label{thm:main}
For every closed smooth manifold \(M\) of even dimension \(2N\ge4\),
there is a \(C^\infty\)-open dense set \(\mathcal U\) in the space
\(\Met^\infty(M)\) of smooth Riemannian metrics such that, for every
\(g\in\mathcal U\) and every \(1\le k<N\),
\begin{equation}\label{eq:main-conclusion}
 \ker(\delta_gQ_{k,g}d)=\ker d.
\end{equation}
\end{theorem}

Write
\begin{align}
 \HH_g^k&=\{e\in\Om^k(M):de=0,\ \delta_gQ_{k,g}e=0\},\notag\\
 J_g^k&:\HH_g^k\longrightarrow H^k_{\mathrm{dR}}(M;\mathbb R),
                 \qquad J_g^k(e)=[e],\label{eq:harmonic-map}\\
 \ker J_g^k&=\HH_g^k\cap d\Om^{k-1}(M).\notag
\end{align}
Strong regularity in degree \(k-1\) is equivalent to \(\ker J_g^k=0\).
The quotient of detour cohomology by the image of de Rham cohomology
is naturally isomorphic to \(\ker L_g/\ker d\), since both groups are
formed modulo the same exact forms \cite[Proposition 2.5]{BG}.

\begin{corollary}[Existence and uniqueness of conformal harmonic representatives]\label{cor:bg-conjecture}
Let \(M\) and \(\mathcal U\) be as in Theorem~\ref{thm:main}.
For every \(g\in\mathcal U\) and every \(1\le k<N\),
\begin{equation}\label{eq:bg-corollary}
 J_g^k:\HH_g^k\xrightarrow{\ \cong\ }H^k_{\mathrm{dR}}(M;\mathbb R).
\end{equation}
\end{corollary}
Aubry and Guillarmou \cite{AG} constructed the operators using the
geometry at infinity. The four-dimensional Maxwell gauge and its
cohomological properties go back to Eastwood and Singer
\cite{ESgauge,ES}. Metric perturbations of zero eigenspaces were used
for the Yamabe operator in \cite{GHJL}. For the ambient formula
for the Hessian of total \(Q\)-curvature at an obstruction-flat metric,
see \cite[(0-4) and Theorem 3.4]{Matsumoto}.

Gover and \v{S}ilhan \cite[Theorem 6.4]{GS} establish surjectivity
on compact conformally Einstein manifolds. They obtain strong
regularity when the Einstein representative is Ricci-flat or has
positive scalar curvature. Aubry and Guillarmou
\cite[Proposition 1.4]{AG} obtain surjectivity from nonnegativity of
\(Q_{k,g}\) on closed forms. {We ask whether small changes in the metric can eliminate
all nonzero exact conformal harmonic forms. Our argument applies
regardless of the sign of the pairing.} Canzani's
generic simplicity results \cite{Canzani} concern nonzero eigenvalues
of conformally covariant elliptic operators. The exterior derivative $d$ in
\(L_g\) forces all closed forms into its kernel at every metric.
Exterior differentiation identifies the remaining kernel
\(\ker L_g/\ker d\) with the exact conformal harmonics.

For \(k=1\), \(\delta_gQ_{1,g}d\) is a nonzero
multiple of the critical scalar GJMS operator \(P_{N,g}\)
\cite{GJMS}, by \cite[Theorems 2.1 and 2.8]{BG}.
Thus Theorem~\ref{thm:main} also gives generic triviality
of the kernel of the critical GJMS operator modulo constants.

This conclusion can be formulated directly in terms of
conformal structures. Let \(\mathcal C(M)\) denote the space
of smooth conformal classes on \(M\), equipped with the
quotient topology induced by the \(C^\infty\) topology on
\(\Met^\infty(M)\). Set
\[
\mathcal T_M
=\mathcal C(M)/\operatorname{Diff}_0(M),
\]
where \(\operatorname{Diff}_0(M)\) is the group of
diffeomorphisms isotopic to the identity.
Strong regularity is invariant under conformal rescaling
and diffeomorphism pullback. Consequently, the genericity
statement of Theorem~\ref{thm:main} descends to these
quotient spaces.

\begin{corollary}\label{cor:conformal-classes}
The set of conformal classes represented by metrics that
are strongly regular in every degree \(k-1\), \(1\le k<N\),
is open and dense in \(\mathcal C(M)\). Its image in
\(\mathcal T_M\) is also open and dense.
In particular, if \(M\) is connected, the critical scalar
GJMS operator \(P_{N,g}\) of order \(2N\) satisfies
\(\ker P_{N,g}=\mathbb R\) on an open dense subset of
\(\mathcal T_M\).
\end{corollary}

The last assertion settles Conjecture~G of
Canzani, Gover, Jakobson, and Ponge \cite{CGJP}.
In dimension four, the critical GJMS operator is the
Paneitz operator \cite{Paneitz},
\[
P_{2,g}u
=\Delta_g^2u+
\delta_g\left(
\tfrac23\operatorname{Scal}_g\,du
-2\operatorname{Ric}_g(\nabla_g u,\cdot)
\right).
\]
Singer \cite[Section~3, Proposition~5 and the discussion
that follows]{Singer1995} proved openness of strong
regularity in dimension four and anticipated its genericity.

Singer also exhibited metrics for which strong regularity
fails. In his example \cite[Section~4]{Singer1995},
\(\Sigma\) is a closed hyperbolic surface and
\(S^2\times\Sigma\) carries the product metric with
curvatures \(+1\) and \(-1\). On functions pulled back
from \(S^2\), the Paneitz operator reduces to
\(\Delta_{S^2}(\Delta_{S^2}-2)\).
The degree-one spherical harmonics therefore yield
nonconstant kernel elements. Nevertheless, by twisting
the sphere factor using irreducible \(SU(2)\)
representations of \(\pi_1(\Sigma)\), Singer constructed
strongly regular, locally conformally flat structures
arbitrarily close to the product structure.

\subsection*{Proof strategy}
In Section~\ref{sec:elliptic}, we use the Hodge decomposition for a
fixed reference metric to construct a self-adjoint elliptic family
on a fixed Sobolev domain. Exterior differentiation identifies its
kernel with \(\ker J_g^k\). Write
\(\nu_k(g)=\dim\ker J_g^k\). If \(\nu_k\) is locally constant, one can
continue a kernel element \(\alpha\) smoothly along any metric variation.
Differentiating the kernel equation twice shows that a reduced Hessian
vanishes in every metric direction. It is the Hessian of
\((Q_{k,g}e,e)_g\) with \(e=d\alpha\) fixed, minus the Green term that
accounts for the motion of \(\alpha\).

Lemma~\ref{lem:leading-variation-coefficients} reduces the principal-order
calculation to a flat metric with \(e\) constant.
Section~\ref{sec:formal} states the variation formulas for a Fourier mode
on a flat torus, and Section~\ref{sec:periodic-calculations} proves them.
The Green correction contributes at the same order as the fixed-form
Hessian and cancels the contribution with coefficient \((N/k)^2\)
(Proposition~\ref{prop:reduced-hessian}).
At each point where \(e\ne0\), one can choose a frequency and a transverse,
trace-free variation for which the resulting symbol is nonzero
(Lemma~\ref{lem:polarization}). This contradicts the vanishing of the
reduced Hessian at a metric with positive, locally constant nullity.

In any nonempty open set of metrics, choose \(g_0\) where \(\nu_k\) is
smallest. Upper semicontinuity gives \(\nu_k(g)<\nu_k(g_0)+1\) near
\(g_0\). Since \(\nu_k\) takes integer values and \(\nu_k(g_0)\) is
the minimum in the open set, \(\nu_k\) is constant near \(g_0\).
The preceding argument yields \(\nu_k(g_0)=0\), proving density.
Upper semicontinuity gives openness, and the intersection over
the finitely many degrees gives Theorem~\ref{thm:main}.

\section*{AI usage statement}

The authors used ChatGPT 6.1 Sol to assist with calculations for the
periodic model in Sections~\ref{sec:formal} 
and to explore proof strategies for some estimates in Section~\ref{sec:elliptic} and \ref{sec:periodic-calculations}. The authors independently
verified and take full responsibility for the
correctness and content of the manuscript.

\section{Preliminaries}\label{sec:preliminaries}

\subsection{The detour complex and notation}

We fix a degree \(1\le k<N\) when working with a single operator.
As described before (7) in Section~2 of \cite{BG}, the weighted dual
of the bundle of \(r\)-forms is \(\Lambda^rT^*M\otimes\mathcal E[2r-2N]\),
where \(\mathcal E[w]\) is the density bundle of conformal weight \(w\).
A representative metric \(g\) trivializes \(\mathcal E[w]\). In this
trivialization the integrated duality pairing is integration of
the ordinary pointwise inner product against \(dV_g\).

In the scale \(g\), the Branson--Gover detour complex in degree \(k-1\)
takes the form
\cite[Theorem~2.1]{BG} and \cite[Theorem~2.8]{BG}
\[
 0\longrightarrow\Om^0\xrightarrow{d}\cdots\xrightarrow{d}\Om^{k-1}
 \xrightarrow{L_g}\Om^{k-1}\xrightarrow{\delta_g}\cdots
 \xrightarrow{\delta_g}\Om^0\longrightarrow 0 .
\]
The right half lists the ordinary \(g\)-representatives of sections
of the weighted dual bundles. Thus, for example, \(Q_{k,\gamma}\)
and \(L_\gamma\) use the scale \(\gamma\) for their outputs.
For \(k=1\), the complex reduces to \(L_g\) between two spaces of
functions. The cohomology at the left
copy of \(\Om^{k-1}\) is \(\ker L_g/d\Om^{k-2}\) for \(k\ge2\).
For \(k=1\) it is \(\ker L_g\).

Let \(dV_g\) denote the Riemannian volume density of \(g\), and
write \(\langle u,v\rangle_{g,x}\) and \((u,v)_g\) for the pointwise
and integrated pairings. We use the nonnegative Hodge Laplacian
\(\Delta_g=d\delta_g+\delta_gd\). For the spectral and symbol
calculations we complexify the bundles.
For a covector \(\xi\), its metric dual \(\xi^{\sharp_g}\) satisfies
\(g(\xi^{\sharp_g},\cdot)=\xi\). We write \(\xi^\sharp\) when the
metric is clear, and \(\iota_X\) for contraction by a vector
\(X\). For an operator \(D\) of order at most \(j\),
\(\sigma_j(D)(x,\zeta)\) denotes its homogeneous symbol of degree \(j\),
for \(\zeta\in T_x^*M\setminus\{0\}\), with the convention
\(\partial_{x_a}\mapsto i\zeta_a\).

For an integer \(s\ge0\), let \(H^s\Om^r\) be the completion of
\(\Om^r\) in the norm
\[
 \|u\|_{H^s,g}^2=\sum_{j=0}^s\int_M |\nabla_g^j u|_g^2\,dV_g ,
\]
where \(\nabla_g\) is the Levi-Civita connection. In particular,
\(L^2_g\Om^r=H^0\Om^r\) with the \(L^2_g\) pairing. On the closed
manifold \(M\), smooth metrics give equivalent Sobolev norms, so
we omit \(g\) from the notation \(H^s\Om^r\).
When the form degree is clear from context, we write \(H^s\) and
\(L_g^2\) for \(H^s\Om^r\) and \(L_g^2\Om^r\), respectively.

\subsection{Ambient construction and normalization}\label{subsec:canonical-definition}
We use the ambient construction of \(Q_{k,g}\) to compute its metric
variations. The construction extends a closed form off the boundary,
applies a power of the ambient Hodge Laplacian, and extracts a boundary
component.

Fefferman and Graham introduced the ambient metric in \cite{FG1985}.
By \cite[Theorem~2.9]{FGambient} and \cite[Theorem~4.5]{FGambient},
every boundary metric has an even formal Poincar\'e metric in normal form,
\[
 g_+=r^{-2}(dr^2+g_r),\qquad g_r|_{r=0}=g,\qquad
 \operatorname{Ric}(g_+)+2Ng_+=O(r^{2N-2}).
\]
The coefficients of \(g_r\) below order \(r^{2N}\) are uniquely determined.
The computations use the jet through \(r^{2N-2k}\), with
\(2(N-k)\le 2N-2\).

Use \(\rho=r^2/2\), and write \(g_{\sqrt{2\rho}}\) for the even
family \(g_r\) expressed in this variable. With the conventions used
before (6.1) in Section~6 of \cite{AG}, the cone metric
\(-d(rt)^2+(rt)^2g_+\) becomes
\[
 \widetilde g=-2\rho\,dt^2-2t\,dt\,d\rho+t^2g_{\sqrt{2\rho}}.
\]
The Fefferman--Graham normal coordinate in equation (3.14) of
\cite{FGambient} is \(-\rho\) in this convention.
Let \(T=t\partial_t\) and \(Y=-t^{-2}\partial_\rho\).
Write \(\widetilde\delta\) and \(\widetilde\nabla\) for the
codifferential and Levi-Civita connection of \(\widetilde g\).
Its Hodge Laplacian is
\(\widetilde\Delta=d\widetilde\delta+\widetilde\delta d\).
For a closed \(k\)-form \(e\), use its coordinate pullback to the
ambient space and set
\[
 F_0=t\,dt\wedge d\rho\wedge e,\qquad
 F_s=\widetilde\Delta^sF_0.
\]
The ambient construction uses a degree-zero homogeneous extension
of the given ordinary form on the cone. The coordinate pullback has
\(\mathcal L_Te=0\), \(\iota_Te=\iota_Ye=0\), and the required
restriction at \(t=1,\rho=0\). In these coordinates
\(\widetilde g(T,\cdot)=\tfrac12 d(\widetilde g(T,T))\) and
\(\mathcal L_T\widetilde g=2\widetilde g\), so
\(\widetilde\nabla T=\Id\).
On \(k\)-forms this gives \(\mathcal L_T=\widetilde\nabla_T+k\), hence
\(\widetilde\nabla_Te=-ke\), as in the paragraph preceding (6.1)
in Section~6 of \cite{AG}.
Since \(F_0\) is closed, commutation with \(d\) gives
\(F_{s+1}=d\widetilde\delta F_s\). The initial form has dilation
degree \(2\), and each ambient Laplacian lowers that degree by two.
Thus \(F_s\) has form degree \(k+2\) and dilation degree \(2-2s\),
with \(\widetilde\nabla_TF_s=-(k+2s)F_s\).
In the ambient formula in Section~6 of \cite{AG}, the initial wedge
factor \(\tfrac12 d(\widetilde g(T,T))\wedge dt/t\) is
\(t\,dt\wedge d\rho\). Since \(Y\) is proportional to
\(\partial_\rho\), the contraction \(\iota_Y\iota_TF_{N-k}\) is
annihilated by \(\iota_{\partial_\rho}\). Dividing the formula
after (6.2) of \cite{AG} by \(2(N-k+1)\) gives
\begin{equation}\label{eq:ambient-definition}
 \begin{aligned}
 Q_{k,g}e&=-\left.\iota_Y\iota_TF_{N-k}\right|_{t=1,\,\rho=0},\\
 Q_k^{\rm BG}e&=2(N-k+1)Q_{k,g}e.
 \end{aligned}
\end{equation}
We use the first line of \eqref{eq:ambient-definition} throughout.
The second uses the normalization of \cite[equation (45)]{BG}.
In this normalization, \cite[Theorem~2.8(iii)]{BG} gives
\begin{equation}\label{eq:bg-detour-normalization}
 L^{\rm BG}_{k-1}=(2N-2k+4)\delta_g Q_k^{\rm BG}d
                 =4(N-k+1)(N-k+2)L_g.
\end{equation}
The scalar factor leaves all the kernels in this paper unchanged.
At \(t=1,\rho=0\), the two contractions and the sign in
\eqref{eq:ambient-definition} extract the coefficient of \(dt\wedge
d\rho\) in \(F_{N-k}\), with \(dt\) and \(d\rho\) placed first.
For the flat normalization check, use orthonormal boundary
coordinates \((x,y_1,\ldots,y_{2N-1})\). Take a constant \((k-1)\)-form
\(\beta\) with \(\iota_{\partial_x}\beta=0\), and set
\(\Theta=dt\wedge d\rho\wedge dx\wedge\beta\).
At the flat metric, raise the \(t,\rho,x\) indices and take the
divergence using the density \(t^{2N+3-2k}\) for this component. This gives
\[
 \widetilde\delta\bigl(t^{\omega-1}f(x)\Theta\bigr)
   =-t^{\omega-3}
    \bigl(f'(x)\,dt\wedge d\rho+(2N-2k+\omega-2)f(x)\,dt\wedge dx\bigr)
       \wedge\beta .
\]
The input form is closed, so applying \(d\) gives
\[
 \widetilde\Delta\bigl(t^{\omega-1}f(x)\Theta\bigr)
       =-t^{\omega-3}f''(x)\Theta .
\]
Starting with \(F_0=tf(x)\Theta\) and iterating the flat formula gives
\[
 Q_{k,g_0}\bigl(f(x)\,dx\wedge\beta\bigr)
   =\bigl((-\partial_x^2)^{N-k}f(x)\bigr)\,dx\wedge\beta ,
\]
where \(g_0\) is the flat boundary metric. Branson and Gover identify
\((d\delta_g)^{N-k}\) as the principal part of \(Q_{k,g}\) on closed
forms up to a nonzero constant in equation (3) of \cite{BG}.
The flat calculation
fixes this constant to be one, so \(L_g\) has principal part
\((\delta_gd)^{N-k+1}\).

The ambient formula defines a natural differential operator whose
coefficients depend smoothly on finitely many jets of \(g\). It uses
only finitely many coefficients of \(g_r\), each determined by the
normal-form equations from the boundary metric
\cite[Theorem~2.9]{FGambient}, \cite[Theorem~4.5]{FGambient}.
See also the universal-polynomial description preceding
\cite[Theorem~2.8]{BG}.

\section{Elliptic realization and constant nullity}\label{sec:elliptic}

Put \(m=N-k\) and fix a reference metric \(g\). Hodge decomposition
gives the orthogonal splitting \(L^2_g\Om^{k-1}=C_g\oplus\ker d\), where \(C_g\) is the
\(L^2_g\)-closure of \(\delta_g\Om^k(M)\); the same notation denotes the
corresponding spaces of smooth forms when appropriate.
For a nearby metric \(\gamma\), write \((\gamma^{-1}g)X=(g(X,\cdot))^{\sharp_\gamma}\).
The positive multiplication operator on \(r\)-forms
\[
 (W_r(\gamma)u)(X_1,\ldots,X_r)
 =(dV_\gamma/dV_g)\,
 u\bigl((\gamma^{-1}g)X_1,\ldots,(\gamma^{-1}g)X_r\bigr)
\]
satisfies \((u,v)_\gamma=(W_r(\gamma)u,v)_g\). Hence
\begin{equation}\label{eq:fixed-adjoint}
 W_{k-1}(\gamma)\delta_\gamma=\delta_gW_k(\gamma).
\end{equation}
Indeed, pairing either side applied to a \(k\)-form \(u\) with a
\((k-1)\)-form \(v\) gives \((u,dv)_\gamma\).

For closed \(u,v\), set \(q_\gamma(u,v)=(Q_{k,\gamma}u,v)_\gamma\).
Let \(\Pi_g^r\) be the orthogonal projection onto
\(\ker\Delta_g|_{\Om^r}\). We add a gauge on the closed summand of the
fixed splitting, as in \cite[proof of Proposition 1.4]{AG}, and set
\begin{equation}\label{eq:elliptic-family}
 \mathcal A_\gamma=\delta_gW_k(\gamma)Q_{k,\gamma}d
                +(d\delta_g)^{m+1}+\Pi_g^{k-1} .
\end{equation}

\begin{lemma}\label{lem:realization}
For every \(\gamma\) in a \(C^\infty\) neighborhood of \(g\),
the operator \eqref{eq:elliptic-family} is self-adjoint and elliptic
with common domain \(H^{2m+2}\Om^{k-1}\). It is block diagonal for
\(C_g\oplus\ker d\), with invertible closed block
\(\Delta_g^{m+1}+\Pi_g^{k-1}\). Exterior differentiation induces an isomorphism
\[
 d:\ker\mathcal A_\gamma\xrightarrow{\ \cong\ }\ker J_\gamma^k.
\]
At the reference metric,
\(\sigma_{2m+2}(\mathcal A_g)(x,\xi)=|\xi|_g^{2m+2}\Id\).
\end{lemma}

\begin{proof}
The first summand of \(\mathcal A_\gamma\) takes values in \(C_g\)
and vanishes on \(\ker d\). Its pairing with \(u,v\) is
\(q_\gamma(du,dv)\), which is symmetric. The remaining summands
vanish on \(C_g\) and act on \(\ker d\) as
\(\Delta_g^{m+1}+\Pi_g^{k-1}\). This block is symmetric and invertible,
since it is the identity on harmonic forms and acts by
\(\lambda^{m+1}\) on each positive eigenspace of \(\Delta_g\).
Thus \(\mathcal A_\gamma\) is block diagonal and formally symmetric.

At \(\gamma=g\), the principal part of the first summand is
\((\delta_gd)^{m+1}\). Adding \((d\delta_g)^{m+1}\) gives
\(\Delta_g^{m+1}\), since \(d^2=\delta_g^2=0\), and
\(\Pi_g^{k-1}\) is smoothing. Hence
\[
 \sigma_{2m+2}(\mathcal A_g)(x,\xi)=|\xi|_g^{2m+2}\Id.
\]
Continuity of the principal symbol and compactness of the unit
cosphere bundle give ellipticity for \(\gamma\) near \(g\).
Elliptic regularity gives the maximal \(L^2_g\) domain
\(H^{2m+2}\Om^{k-1}\), with graph norm equivalent to the Sobolev norm.
Smooth forms are a core, so formal symmetry gives self-adjointness
on that domain. The Sobolev embedding gives compact resolvent,
and the kernel consists of smooth forms.

On \(C_g\), \eqref{eq:fixed-adjoint} gives
\(\mathcal A_\gamma=W_{k-1}(\gamma)L_\gamma\).
By Hodge decomposition every smooth exact \(k\)-form \(e\) has a
unique primitive \(\alpha\in C_g\). Since \(W_{k-1}(\gamma)\) and the
closed block are invertible, \(\mathcal A_\gamma\alpha=0\) is
equivalent to \(\delta_\gamma Q_{k,\gamma}e=0\) for this primitive.
This gives the asserted isomorphism. The finite-jet dependence
on \(\gamma\) also gives continuity in the \(C^\infty\) metric
topology and smooth dependence along smooth metric paths in
operator norm from \(H^{2m+2}\Om^{k-1}\) to \(L^2_g\).
\end{proof}

Set \(K_g=\ker\mathcal A_g\), and let \(G\) be the inverse of
\(\mathcal A_g\) on \(K_g^\perp\), extended by zero on \(K_g\).
It is bounded on \(L^2_g\) and preserves \(C_g\). When \(k=1\),
the term \(d\delta_g\) is zero and the closed block consists of the
harmonic projection onto locally constant functions.

\subsection{Metric variations and constant nullity}\label{sec:constant-nullity}

Put \(\nu_k(g)=\dim\ker J_g^k=\dim K_g\). For a smooth real symmetric
tensor \(h\) and \(\alpha\in K_g\), write
\[
 e=d\alpha,\qquad
 \mathcal F_\alpha h=\left.\frac{d}{dt}\right|_0\mathcal A_{g+th}\alpha
     =\left.\frac{d}{dt}\right|_0 L_{g+th}\alpha.
\]
The last equality follows from \eqref{eq:fixed-adjoint},
\(W_{k-1}(g)=\Id\), and \(L_g\alpha=0\).
If the nullity is locally constant, \(\alpha\) extends to a smooth
family \(\alpha_t\in\ker\mathcal A_{g+th}\) with \(\alpha_0=\alpha\).
The proof of Lemma~\ref{lem:constant-nullity} shows that the first variation
with \(e\) fixed vanishes. The \(K_g^\perp\) component of the derivative
of \(\alpha_t\) contributes through \(G\) to the reduced second variation
\begin{equation}\label{eq:second-variation}
 \mathscr S_g(h;\alpha)
 :=\left.\frac{d^2}{dt^2}\right|_0q_{g+th}(e,e)
                         -2(\mathcal F_\alpha h,G\mathcal F_\alpha h)_g.
\end{equation}

\begin{lemma}\label{lem:constant-nullity}
The function \(\nu_k\) is finite and upper semicontinuous. If it is
constant on a neighborhood of \(g\), then
\(\mathscr S_g(h;\alpha)=0\) for every smooth real symmetric tensor
\(h\) and every \(\alpha\in K_g\).
\end{lemma}
\begin{proof}
Choose a small circle about zero containing no other spectral value of
\(\mathcal A_g\). The common domain and operator-norm smoothness from
Lemma~\ref{lem:realization} imply that its resolvent contour defines,
for nearby \(\gamma\), a spectral projection \(E_\gamma\), smooth as
an operator \(L^2_g\to H^{2m+2}\), of rank \(\dim K_g\).
Since \(\ker\mathcal A_\gamma\subset\im E_\gamma\),
this proves upper semicontinuity. If the nullity is constant, then
\(\im E_\gamma=\ker\mathcal A_\gamma\).

For \(g_t=g+th\), write \(\mathcal A_t=\mathcal A_{g_t}\) and set
\(\alpha_t=E_{g_t}\alpha\). Then
\(\alpha_0=\alpha\) and \(\mathcal A_t\alpha_t=0\).
Differentiation at zero gives
\[
 \mathcal A_g\dot\alpha=-\mathcal F_\alpha h,\qquad
 \mathcal F_\alpha h\perp K_g,\qquad
 \dot\alpha+G\mathcal F_\alpha h\in K_g.
\]
Differentiating \(\mathcal A_t\alpha_t=0\) again and pairing with
\(\alpha\) gives
\[
 0=(\mathcal A''_0\alpha,\alpha)_g
       +2(\mathcal F_\alpha h,\dot\alpha)_g
   =\left.\frac{d^2}{dt^2}\right|_0q_{g_t}(e,e)
       -2(\mathcal F_\alpha h,G\mathcal F_\alpha h)_g.
\]
The last equality uses \((\mathcal A_{g_t}\alpha,\alpha)_g
=q_{g_t}(e,e)\) and \(\mathcal F_\alpha h\perp K_g\).
\end{proof}

\subsection{Locality of the principal terms}

To evaluate the local variations entering the reduced Hessian, we
first determine which background data can enter at principal order.

\begin{lemma}\label{lem:leading-variation-coefficients}
For a fixed closed \(k\)-form \(e\), the second metric variation of
\(\langle Q_{k,g}e,e\rangle_{g,x}\,dV_g\) has total differential order
at most \(2m\) in the two metric variations. The first metric
variation of \(\delta_gQ_{k,g}e\) has order at most \(2m+1\).
In local coordinates the coefficients at these orders depend only
on \(g(x)\) and \(e(x)\).
\end{lemma}
\begin{proof}
Choose coordinates centered at \(0\) and write \(\phi_\lambda(x)=\lambda x\).
Set \(g_\lambda=\lambda^{-2}\phi_\lambda^*g\) and
\(e_\lambda=\lambda^{-k}\phi_\lambda^*e\), so their coordinate
components are those of \(g\) and \(e\) evaluated at \(\lambda x\).
Naturality, the constant-scale law for \(Q_{k,g}\), and
\(\delta_{c^2g}=c^{-2}\delta_g\) give
\[
 Q_{k,g_\lambda}e_\lambda=\lambda^{2m}(Q_{k,g}e)(\lambda x),\qquad
 \delta_{g_\lambda}Q_{k,g_\lambda}e_\lambda
                 =\lambda^{2m+1}(\delta_gQ_{k,g}e)(\lambda x).
\]
At \(0\), each order-\(r\) jet scales by \(\lambda^r\). The
coordinate coefficient of the pairing density and the components of
\(\delta_gQ_{k,g}e\) thus have weights \(2m\) and \(2m+1\), respectively.

Choose a common finite jet order \(s\ge1\) for both expressions.
Represent an \(s\)-jet of \(e\) by its Taylor polynomial
\(\sum_{r=0}^s e_{[r]}\), where \(e_{[r]}\) has homogeneous
coefficients of degree \(r\). The condition \(j_0^{s-1}(de)=0\)
is exactly \(de_{[r]}=0\) for \(1\le r\le s\), since exterior
differentiation lowers the coefficient degree by one.
The constant term is closed as well. For the radial field
\(R=\sum_i x^i\partial_{x^i}\), Cartan's formula gives
\[
 d\iota_R e_{[r]}=\mathcal L_R e_{[r]}=(r+k)e_{[r]}.
\]
Thus the prescribed jet is realized by the exterior derivative of
\(\sum_{r=0}^s(r+k)^{-1}\iota_R e_{[r]}\). Cutting off this
primitive inside the chart, with the cutoff equal to one near \(0\),
gives a global exact \(k\)-form with the same jet. The constraints
are independent of \(g\) and define a linear subspace at each weight
\(r\), which is preserved by scaling by \(\lambda^r\).

With \(g(0)\) and \(e(0)\) fixed, let \(z=(z_r)_{r=1}^s\)
denote linear coordinates on the allowed positive-order metric and
form jets. Write \(S_\lambda z=(\lambda^r z_r)_{r=1}^s\).
By Section~\ref{subsec:canonical-definition}, a coordinate component
of weight \(d\) is a smooth function \(F\) of these jets, with
\(F(S_\lambda z)=\lambda^d F(z)\). Taylor expansion at \(z=0\)
through ordinary degree \(d\) gives
\[
 F(z)=\sum_{|\alpha|\le d}c_\alpha z^\alpha+O(|z|^{d+1})
 \qquad (z\to0).
\]
Here \(\alpha_r\) is the part of the multi-index in the block \(z_r\).
Since \(S_\lambda z=O(\lambda)\) for each fixed \(z\), the expansion
at \(S_\lambda z\) gives, as \(\lambda\downarrow0\),
\[
 \lambda^d F(z)=\sum_{w=0}^d\lambda^w P_w(z)+O(\lambda^{d+1}),
 \qquad
 P_w(z)=\sum_{\substack{|\alpha|\le d\\ \sum_r r|\alpha_r|=w}}
             c_\alpha z^\alpha.
\]
The polynomial in \(\lambda\) obtained by subtracting
\(\lambda^d F(z)\) from the sum has degree at most \(d\) and
vanishes to order at least \(d+1\). Hence \(P_w(z)=0\) for \(w<d\)
and \(F(z)=P_d(z)\). The coefficients depend smoothly on
\(g(0)\) and \(e(0)\).
Metric differentiation replaces jets by jets of the variations of
the same weight. Terms with a total of \(d\) derivatives on the
variations exhaust the weight, leaving only \(g(0)\) and \(e(0)\)
in their coefficients.
\end{proof}

\section{The periodic model and variation formulas}\label{sec:formal}
By the jet weights in Lemma~\ref{lem:leading-variation-coefficients},
terms with fewer derivatives on the metric variations contain a
positive-order jet of the background metric or of \(e\). In flat
coordinates with \(e\) constant, these jets vanish, leaving just
the terms of principal order.
A single Fourier mode tests the principal symbol in its cotangent
direction \(\xi\).
For a unit covector \(\xi\), symmetric endomorphisms \(H\) with
\(H\xi^\sharp=0\) and \(\operatorname{tr}H=0\) give a single radial
profile \eqref{eq:linear-profile}. For each nonzero \(e\),
Lemma~\ref{lem:polarization} gives a pair \((\xi,H)\) in this family
for which the principal symbol of the reduced Hessian is nonzero.

For a symmetric endomorphism \(H\), an \(r\)-form \(\varphi\) with
\(r\ge1\), and tangent vectors \(v_1,\ldots,v_r\), define
\[
 (H^\sharp\varphi)(v_1,\ldots,v_r)
 =\sum_{\ell=1}^r\varphi(v_1,\ldots,Hv_\ell,\ldots,v_r).
\]
Set \(H^\sharp f=0\) for functions \(f\). Endomorphism norms are
Hilbert--Schmidt norms.

The variation formulas use the constants
\[
 \chi=\frac Nk,\qquad a=\frac{k}{2(2k+1)},\qquad
 B=\frac{m!(2k)_{m}}{2(k)_{m}^2}>0,
\]
where \((z)_q=z(z+1)\cdots(z+q-1)\) and \((z)_0=1\).

\subsection{The formal Einstein metric}\label{subsec:formal-einstein}

For the periodic calculation, we work on the flat torus
\(\mathbb T^{2N}\), with periodic coordinates
\(x,y_1,\ldots,y_{2N-1}\). Fix a diagonal trace-free matrix
\[
H=\operatorname{diag}(h_1,\ldots,h_{2N-1}),
\qquad
\sum_{a=1}^{2N-1}h_a=0,
\qquad
w=\operatorname{tr}(H^2)
=\sum_{a=1}^{2N-1}h_a^2.
\]
Consider the one-parameter family of boundary metrics
\[
g_\varepsilon
=dx^2+\sum_{a=1}^{2N-1}
(1+\varepsilon h_a\cos x)\,dy_a^2,
\]
for sufficiently small \(|\varepsilon|\).
This is a single Fourier-mode, trace-free perturbation
of the flat metric.

Throughout this calculation, we write
\[
f^{(q)}=[\varepsilon^q]f
\]
for the coefficient of \(\varepsilon^q\) in the Taylor
expansion of \(f\). In particular,
\[
f^{(2)}
=\frac12\left.
\partial_\varepsilon^2 f
\right|_{\varepsilon=0}.
\]
Unsubscripted angle brackets denote averaging
over the periodic \(x\)-variable:
\[
\langle f\rangle
=\frac{1}{2\pi}\int_0^{2\pi}f(x)\,dx.
\]

The formal geodesic metric has
\(g_r=u\,dx^2+\sum v_a\,dy_a^2\). Its isometries force this diagonal
ansatz below the undetermined order. Linearizing the Einstein equations
as in Section~\ref{app:coordinate-einstein} gives
\begin{align}
 &u^{(1)}=0,\qquad v_a^{(1)}=h_a A(r)\cos x, \notag\\
 \shortintertext{where the radial profile \(A\) satisfies}
 &A''-\frac{2N-1}{r}A'-A=O(r^{2m}). \notag\\
 \shortintertext{With \(A(0)=1\), matching powers of \(r\) gives}
 &[r^{2q}]A=\frac{(-1)^q\Gamma(N-q)}{4^q q!\Gamma(N)},
                       \qquad 0\le q\le m. \label{eq:linear-profile}
\end{align}
Put \(\lambda_a=\tfrac12\log v_a\), \(\Lambda=\sum\lambda_a\).
Write the second means as
\[
 \ell_a=\langle\lambda_a^{(2)}\rangle,\quad
 U=\langle u^{(2)}\rangle.
\]
Their boundary data are \(\ell_a(0)=-h_a^2/8\), \(U(0)=0\).
Put \(\mathfrak F=A^2-(A')^2\) and
\(\mathfrak F_q=[r^{2q}]\mathfrak F\).

\subsection{The periodic variation formulas}\label{sec:ambient-variation}

Let \(g_\varepsilon\) be the periodic boundary metric of
Section~\ref{sec:formal}, with \(g_0=dx^2+\sum_a dy_a^2\).
For a constant closed \(k\)-form \(e\), define the unreduced Hessian by
\[
 R(e;dx,H)=\left.\frac{d^2}{d\varepsilon^2}\right|_0
 \frac{1}{(2\pi)^{2N}}\int_{\mathbb T^{2N}}
     \langle Q_{k,g_\varepsilon}e,e\rangle_{g_\varepsilon,x}\,dV_{g_\varepsilon}.
\]
This normalized integral is the period mean used below.

\begin{proposition}\label{prop:raw-hessian}
Let \(e=dx\wedge\beta+\eta\) be constant, with both components
transverse to \(dx\), and let \(H\) be diagonal in the model
coordinates, transverse to \(dx\), and trace-free. The first forcing is
\begin{equation}\label{eq:firstforcing}
 \left.\frac{d}{d\varepsilon}\right|_0
       \bigl(\delta_{g_\varepsilon}Q_{k,g_\varepsilon}e\bigr)
                   =-\chi\sin x\,H^\sharp\beta.
\end{equation}
The unreduced Hessian is
\begin{equation}\label{eq:raw}
 \begin{aligned}
 R(e;dx,H)={}&B|H^\sharp\eta|^2
   +(\chi^2-\chi B)|H^\sharp\beta|^2\\
 &-aB|H|^2|\eta|^2+\chi aB|H|^2|\beta|^2.
 \end{aligned}
\end{equation}
\end{proposition}

For a unit covector \(\xi\) and a transverse trace-free \(H\),
choose an orthonormal coframe with first covector \(\xi\) whose
remaining dual vectors diagonalize \(H|_{\ker\xi}\). Use it as the
coordinate coframe of the model torus and denote the resulting Hessian
by \(R(e;\xi,H)\). In \eqref{eq:raw} the components are then
\(\beta=\iota_{\xi^\sharp}e\) and \(\eta=e-\xi\wedge\beta\).
The norms in \eqref{eq:raw} are independent of this coframe, and
the form \(H^\sharp\beta\) in \eqref{eq:firstforcing} transforms
with the coframe.

We call the \(\eta\) contribution tangential and the
\(dx\wedge\beta\) contribution normal, relative to the \(x\)-direction.
After factoring out the transverse exterior form, the remaining
forms in \((t,\rho,x)\) have degrees two and three, respectively.
Section~\ref{app:variation} proves the formulas, and
Section~\ref{sec:variation-checks} gives two independent checks.

\section{The reduced metric Hessian and generic strong regularity}\label{sec:localization}

Fix the background metric \(g\), and put \(\mathcal V=\Sym^2T^*M\),
with its induced metric and \(L^2_g\) pairing. A symmetric endomorphism
\(H\) represents the covariant tensor \(h_H=g(H\cdot,\cdot)\).
All adjoints below use the fixed \(g\)-pairings. We use the classical
pseudodifferential calculus \(\Psi_{\mathrm{cl}}^r\) on the
complexified bundles, with symbol convention \(\partial_x\mapsto i\zeta\).

Let \(x\in M\), let \(e\in\Lambda^kT_x^*M\), and let
\(\xi\) be a unit covector at \(x\). Set
\(\beta=\iota_{\xi^\sharp}e\) and \(\eta=e-\xi\wedge\beta\).
The interior-product rule, together with \(|\xi|_g=1\) and
\(\iota_{\xi^\sharp}^{\,2}=0\), gives
\[
 e=\xi\wedge\beta+\eta,\qquad
 \iota_{\xi^\sharp}\beta=\iota_{\xi^\sharp}\eta=0.
\]
We consider a symmetric endomorphism \(H\) satisfying
\(H\xi^\sharp=0\) and \(\operatorname{tr}H=0\).

The quadratic expression for the reduced symbol is
\begin{equation}\label{eq:effective-symbol}
 \mathcal J(e;\xi,H)=B\left\{
 |H^\sharp\eta|^2-\chi|H^\sharp\beta|^2
 -a|H|^2\bigl(|\eta|^2-\chi|\beta|^2\bigr)\right\}.
\end{equation}

\subsection{Principal symbols of metric variations}

For \(e\) fixed, let primes on \(W_k\) and \(Q_{k,g}\) denote
derivatives of \(W_k(g_t)\) and \(Q_{k,g_t}\) at \(t=0\),
where \(g_t=g+th\). The product rule gives
\[
 \left.\frac{d^2}{dt^2}\right|_0q_{g_t}(e,e)
   =\bigl(W_k''Q_{k,g}e+2W_k'Q_{k,g}'e+Q_{k,g}''e,e\bigr)_g.
\]
The term \(W_k''Q_{k,g}e\) is algebraic in \(h\). The mixed term
\(2W_k'Q_{k,g}'e\) can enter at principal order.

The Hessian of \(q_\gamma(e,e)\) defines by
integration by parts a real formally self-adjoint operator
\(\mathbb H_{g,e}\) on \(\mathcal V\), of differential order at most
\(2m\), satisfying
\[
 \left.\frac{\partial^2}{\partial s\,\partial t}\right|_{s=t=0}
 q_{g+sh_1+th_2}(e,e)=(h_1,\mathbb H_{g,e}h_2)_g.
\]
In local coordinates a Hessian term
\(C_{A,B}\partial^Ah_1\partial^Bh_2\), with
\(|A|+|B|=2m\), has principal symbol factor
\((-i)^{|A|}i^{|B|}=(-1)^{|A|+m}\).
Evaluating both metric arguments on unit-frequency cosines and averaging
over a period gives half this factor, since
\(\cos((|A|-|B|)\pi/2)=(-1)^{|A|+m}\).

Freeze \(g,e\) at \(x\) as in Lemma~\ref{lem:leading-variation-coefficients}.
For \(\zeta\ne0\), put \(\xi=\zeta/|\zeta|_g\) and take \(H\)
transverse and trace-free. The unit-frequency test of
Section~\ref{sec:ambient-variation}, followed by homogeneity in \(\zeta\),
gives
\begin{equation}\label{eq:raw-metric-symbol}
 \langle h_H,\sigma_{2m}(\mathbb H_{g,e})(x,\zeta)h_H\rangle_{g,x}
           =2|\zeta|_g^{2m}R(e(x);\xi,H).
\end{equation}
The factor two accounts for the cosine-square average.

Fix \(\alpha\in K_g\) and \(e=d\alpha\). The forcing map defined in
Section~\ref{sec:constant-nullity} is a real differential operator
\[
 \mathcal F_\alpha:C^\infty(\mathcal V)\longrightarrow\Om^{k-1}(M).
\]
It takes values in \(C_g\). Since \(\delta_gQ_{k,g}e=0\),
by \eqref{eq:fixed-adjoint}, \(\mathcal F_\alpha h\) is the
derivative of \(\delta_{g_t}Q_{k,g_t}e\) at \(t=0\). Its order is at most \(2m+1\) by
Lemma~\ref{lem:leading-variation-coefficients}, so
\eqref{eq:firstforcing} gives, with \(\beta=\iota_{\xi^\sharp}e(x)\)
and the same \(H\),
\begin{equation}\label{eq:forcing-symbol}
 \sigma_{2m+1}(\mathcal F_\alpha)(x,\zeta)h_H
     =i\chi|\zeta|_g^{2m+1}H^\sharp\beta.
\end{equation}
The factor \(i\) gives
the \(-\sin x\) output when acting on \(\cos x\).
The term \(\delta_g(W_k'Q_{k,g}e)\) has order at most one in the
metric variation and does not enter \eqref{eq:forcing-symbol}.

\subsection{The reduced Hessian}\label{subsec:reduced-hessian}

Section~\ref{sec:constant-nullity} shows that, when the nullity is
locally constant, the component of \(\dot\alpha\) orthogonal to \(K_g\)
is \(-G\mathcal F_\alpha h\). We now compute how this motion contributes
to the principal symbol of the second variation.

On a closed manifold, let \(A\) be a classical elliptic
pseudodifferential operator of order \(r>0\) on a vector bundle \(E\).
If \(A:H^r(E)\to L^2(E)\) is invertible, the elliptic inverse theorem
gives
\[
 A^{-1}\in\Psi_{\mathrm{cl}}^{-r},\qquad
 \sigma_{-r}(A^{-1})(x,\zeta)=\sigma_r(A)(x,\zeta)^{-1}
 \quad(\zeta\ne0).
\]
This is the invertible case of \cite[Chapter 9, Proposition 18]{Melrose}.

At the reference metric \(g\), the spectral projection \(E_g\) from
Lemma~\ref{lem:constant-nullity} is the \(L^2_g\)-orthogonal projection
onto \(K_g\).
Its Schwartz kernel is smooth because \(K_g\) is finite-dimensional
and consists of smooth forms. By Lemma~\ref{lem:realization},
\(\mathcal A_g\) is self-adjoint with compact resolvent.
The operator \(\mathcal A_g+E_g\) acts as the identity on \(K_g\)
and as \(\mathcal A_g\) on \(K_g^\perp\), where zero is separated
from the spectrum. It is therefore invertible from
\(H^{2m+2}\Om^{k-1}\) to \(L^2_g\Om^{k-1}\). It has the same
principal symbol as \(\mathcal A_g\), since \(E_g\) is smoothing.
Subtracting \(E_g\) from its inverse gives the Green operator
of Section~\ref{sec:constant-nullity}, which vanishes on \(K_g\)
and agrees with the inverse of \(\mathcal A_g\) on \(K_g^\perp\):
\[
 G=(\mathcal A_g+E_g)^{-1}-E_g
   \in\Psi_{\mathrm{cl}}^{-2m-2},
 \qquad \mathcal A_gG=\Id-E_g.
\]
Both \((\mathcal A_g+E_g)^{-1}\) and \(E_g\) are real and
self-adjoint, so \(G\) is as well. The right-hand side of
\(\mathcal A_gG=\Id-E_g\) has principal symbol \(\Id\).
Taking principal
symbols and using Lemma~\ref{lem:realization} gives, for \(\zeta\ne0\),
\begin{equation}\label{eq:Green-symbol}
 \sigma_{-2m-2}(G)(x,\zeta)
   =\sigma_{2m+2}(\mathcal A_g)(x,\zeta)^{-1}
   =|\zeta|_g^{-2m-2}\Id.
\end{equation}
Thus the principal symbol of the global Green operator is determined
pointwise by that of \(\mathcal A_g\). The Green term
\(\mathcal F_\alpha^*G\mathcal F_\alpha\) has order at most
\(2(2m+1)-(2m+2)=2m\), so it can contribute at the
same order as \(\mathbb H_{g,e}\).

\begin{proposition}\label{prop:reduced-hessian}
For \(\alpha\in K_g\) and \(e=d\alpha\), the operator
\begin{equation}\label{eq:reduced-hessian}
 \mathcal P_\alpha=\mathbb H_{g,e}
             -2\mathcal F_\alpha^*G\mathcal F_\alpha
       \in\Psi_{\mathrm{cl}}^{2m}(M;\mathcal V)
\end{equation}
is real and formally self-adjoint. For every smooth real metric
variation \(h\),
\begin{equation}\label{eq:reduced-quadratic-form}
 \mathscr S_g(h;\alpha)=(h,\mathcal P_\alpha h)_g.
\end{equation}
For \(\zeta\ne0\), put \(\xi=\zeta/|\zeta|_g\). If \(H\) is
symmetric, trace-free, and kills \(\xi^\sharp\), then
\begin{equation}\label{eq:reduced-symbol}
 \langle h_H,\sigma_{2m}(\mathcal P_\alpha)(x,\zeta)h_H\rangle_{g,x}
       =2|\zeta|_g^{2m}\mathcal J(e(x);\xi,H).
\end{equation}
\end{proposition}
\begin{proof}
All operators act on smooth sections, and the adjoint in
\eqref{eq:reduced-hessian} is taken in the fixed \(g\)-pairings.
The two terms are real and formally self-adjoint, and
\eqref{eq:reduced-quadratic-form} follows from \eqref{eq:second-variation}.
Under complexification, the principal symbol of the formal adjoint
is the Hermitian adjoint. By symbol composition
\cite[Chapters 6--7]{Melrose}, the order-\(2m\) contribution of
the Green term to the quadratic form at \(h_H\) is
\[
 \begin{aligned}
 &\big\langle h_H,
   \sigma_{2m}(2\mathcal F_\alpha^*G\mathcal F_\alpha)(x,\zeta)
      h_H\big\rangle_{g,x}\\
 &\qquad=2|\zeta|_g^{-2m-2}
    \big|i\chi|\zeta|_g^{2m+1}H^\sharp\beta\big|^2\\
 &\qquad=2|\zeta|_g^{2m}\chi^2|H^\sharp\beta|^2,
 \end{aligned}
\]
where the second line uses \eqref{eq:forcing-symbol} and
\eqref{eq:Green-symbol}. The exponent in the last line shows how
the two forcing factors combine with the negative order of \(G\)
to give the same order as the unreduced Hessian. Subtracting the
Green contribution from \eqref{eq:raw-metric-symbol} gives
\[
 \langle h_H,\sigma_{2m}(\mathcal P_\alpha)(x,\zeta)h_H\rangle_{g,x}
 =2|\zeta|_g^{2m}
       \{R(e(x);\xi,H)-\chi^2|H^\sharp\beta|^2\}.
\]
In \eqref{eq:raw}, the coefficient of \(|H^\sharp\beta|^2\)
therefore changes from \(\chi^2-\chi B\) to \(-\chi B\).
This gives \eqref{eq:reduced-symbol} by \eqref{eq:effective-symbol}.
\end{proof}

\subsection{Nonvanishing of the effective symbol}\label{sec:nonvanishing}

\begin{lemma}\label{lem:polarization}
For every nonzero \(k\)-form \(e\) on a Euclidean space of dimension
\(2N\), with \(1\le k<N\), there are a unit covector \(\xi\)
and a transverse trace-free symmetric endomorphism \(H\) such that
\(\mathcal J(e;\xi,H)\ne0\).
\end{lemma}
\begin{proof}
Averaging over \(\xi\) and \(H\) gives a pair with nonzero symbol
for \(k\ge3\). In degrees one and two the average vanishes at
\(N=k+6\), so we choose \(\xi\) and \(H\) directly.
For \(k=1\), choose \(\xi\) parallel to \(e\ne0\) and any
nonzero transverse trace-free \(H\). Then \(\eta=0\),
\(H^\sharp\beta=0\), and
\(\mathcal J/B=N|H|^2|e|^2/6>0\).

For \(k=2\), take an orthonormal coframe in which
\(e=\sum_{j=1}^N\lambda_j u^{2j-1}\wedge u^{2j}\).
Choose \(j\) minimizing \(\lambda_j^2\), put \(\xi=u^{2j-1}\),
and let \(H\) have eigenvalues \(+1,-1\) on another two-plane
of this decomposition and vanish elsewhere. Then \(H\) is transverse
and trace-free, and \(H^\sharp\beta=H^\sharp\eta=0\). Since \(N>2\),
\[
 |\eta|^2-\chi|\beta|^2
   =|e|^2-(1+N/2)\lambda_j^2
   \ge\frac{N-2}{2N}|e|^2>0.
\]
Hence \(\mathcal J=-aB|H|^2(|\eta|^2-\chi|\beta|^2)<0\).

For the average, fix a unit covector \(\xi\) and identify
\((\xi^\sharp)^\perp\) with \(\mathbb R^{2N-1}\).
Averaging \(H\) uniformly over the unit sphere in
\(\Sym_0^2(\mathbb R^{2N-1})\) gives
\[
 \mathbb E_H|H^\sharp\omega|^2=\frac{r(2N-1-r)}{(2N-1)(2N-2)}|\omega|^2,
                                  \qquad\omega\in\Lambda^r\mathbb R^{2N-1}.
\]
For an orthonormal basis \(\{H_a\}\) of this symmetric-tensor space,
\[
 \sum_a(H_a)_{ij}(H_a)_{lm}
 =\frac{\delta_{il}\delta_{jm}+\delta_{im}\delta_{jl}}2
                      -\frac{\delta_{ij}\delta_{lm}}{2N-1}.
\]
Contraction on \(\Lambda^r\mathbb R^{2N-1}\) gives
\(r(2N-1-r)(2N+1)/(2(2N-1))\); division by the dimension
\((2N-2)(2N+1)/2\) gives the stated average.
Averaging the resulting expression over the unit covectors \(\xi\),
and using \(\mathbb E|\beta|^2=k|e|^2/(2N)\), gives
\[
 \frac{\mathbb E_{\xi,H}\mathcal J}{B|e|^2}
 =\frac{m\{6k^2+km(9-2m)+2m\}}
       {4(2k+1)(N-1)N(2N-1)}.
\]
If the numerator vanishes, then \(k\mid 2m\). Set
\(n=2m/k\), a positive integer. Dividing the vanishing factor by
\(km\) gives
\[
 kn=9+\frac{12}{n}+\frac{2}{k}.
\]
For \(k\ge3\) and \(n\ge5\), the left side is at least \(15\) and
the right side is less than \(13\). Hence \(1\le n\le4\), so
\(12/n\) is an integer. Since \(0<2/k<1\), the equality is
impossible. The average is therefore nonzero for \(k\ge3\).
\end{proof}

\subsection{Proof of the main theorem and corollaries}\label{sec:genericity}

\begin{samepage}

\begin{proof}[Proof of Theorem~\ref{thm:main}]
First suppose \(M\) is connected and fix \(1\le k<N\).
In a nonempty open set \(\mathcal O\) of metrics,
choose \(g\) at which \(\nu_k\) attains its minimum. By
Lemma~\ref{lem:constant-nullity}, \(\nu_k\le\nu_k(g)\) in some
neighborhood of \(g\) in \(\mathcal O\), so \(\nu_k\) is constant there.
If \(\nu_k(g)>0\), choose \(0\ne\alpha\in K_g\). The same lemma and
\eqref{eq:reduced-quadratic-form} give
\[
 (h,\mathcal P_\alpha h)_g=0
 \qquad\text{for every smooth real }h.
\]
Real polarization gives \(\mathcal P_\alpha=0\), since this operator is
real and formally self-adjoint.
On the other hand, Lemma~\ref{lem:realization} gives \(d\alpha\ne0\).
At a point \(x\) where \(d\alpha(x)\ne0\),
Lemma~\ref{lem:polarization} and Proposition~\ref{prop:reduced-hessian}
give a nonzero principal symbol of \(\mathcal P_\alpha\), a contradiction.
Thus \(\nu_k(g)=0\), proving that \(\mathcal U_k=\{g:\nu_k(g)=0\}\) is dense.
It is open by upper semicontinuity. The finite intersection
\(\mathcal U=\bigcap_{k=1}^{N-1}\mathcal U_k\) is open and dense,
and \(\ker J_g^k=0\) is equivalent to \(\ker L_g=\ker d\).
For disconnected \(M\), take the finite product of these open dense
sets on its connected components.
\end{proof}
\end{samepage}

\begin{proof}[Proof of Corollary~\ref{cor:bg-conjecture}]
Strong regularity implies surjectivity by \cite[Theorem 2.6]{BG};
here the representative can be obtained directly from the inverse.
Fix \(g\in\mathcal U\) and \(1\le k<N\). By
Lemma~\ref{lem:realization}, \(\mathcal A_g\) is invertible and its
inverse preserves \(C_g\). For a closed \(k\)-form \(\omega\), set
\[
 \alpha=-\mathcal A_g^{-1}\delta_gQ_{k,g}\omega\in C_g,
 \qquad \widehat\omega=\omega+d\alpha.
\]
Since \(\mathcal A_g=\delta_gQ_{k,g}d\) on \(C_g\),
\(\delta_gQ_{k,g}\widehat\omega=0\). Thus \(\widehat\omega\in\HH_g^k\)
and \([\widehat\omega]=[\omega]\), proving surjectivity of \(J_g^k\).
Strong regularity gives \(\ker J_g^k=0\), which proves uniqueness.
\end{proof}

\begin{proof}[Proof of Corollary~\ref{cor:conformal-classes}]
Let \(\mathcal U\) be the set of metrics that are strongly regular in
all these degrees. Lemma~\ref{lem:constant-nullity} and
Theorem~\ref{thm:main} imply that \(\mathcal U\) is open and dense.
The space \(\ker J_g^k\) is unchanged by conformal rescaling because
\(\HH_g^k\) is conformally invariant \cite[Theorem 2.8(ii)]{BG}
and exactness is independent of the metric. Thus \(\mathcal U\) is
saturated under conformal rescaling.
Its image under the quotient map is open by the definition of the
quotient topology. It is dense because the quotient map is continuous
and surjective. For \(\phi\in\operatorname{Diff}_0(M)\), naturality
gives \(\phi^*:\ker J_g^k\cong\ker J_{\phi^*g}^k\). The open dense set in
\(\mathcal C(M)\) is thus saturated under \(\operatorname{Diff}_0(M)\),
and the same quotient argument applies. Taking \(k=1\) gives
\(\ker(\delta_gQ_{1,g}d)=\ker d\) on functions. Since
\(P_{N,g}\) is a nonzero multiple of \(\delta_gQ_{1,g}d\),
its kernel consists of constants when \(M\) is connected.
\end{proof}

\section{Calculations for the periodic model}\label{sec:periodic-calculations}

We prove Proposition~\ref{prop:raw-hessian}. We first derive the
coordinate Einstein equations and compute the averaged metric
coefficients entering the second variation.

\subsection{Coordinate calculations for the formal Einstein metric}\label{app:coordinate-einstein}
We use the notation of Section~\ref{subsec:formal-einstein}.
Set
\(g_B=r^{-2}(dr^2+u\,dx^2)\) and
\(\psi_a=\lambda_a-\log r\), so
\(g_+=g_B+\sum e^{2\psi_a}dy_a^2\).
The only extra connection coefficients, with base indices \(i,l\),
are
\[
 \Gamma^{y_a}_{i y_b}=\delta_{ab}\partial_i\psi_a,\qquad
 \Gamma^i_{y_a y_b}=-\delta_{ab}e^{2\psi_a}
                                      (g_B)^{il}\partial_l\psi_a.
\]
Substitution in the definition of Ricci gives
\[
 \operatorname{Ric}_{ij}=\operatorname{Ric}_{B,ij}
    -\sum_a(\nabla_B^2\psi_a+d\psi_a\otimes d\psi_a)_{ij},
\]
\[
 e^{-2\psi_a}\operatorname{Ric}_{y_a y_a}
     =-\operatorname{tr}_{g_B}\nabla_B^2\psi_a
        -\langle d\psi_a,d\textstyle\sum_b\psi_b\rangle_B,
\]
Writing \(\overline K=-u_{rr}/(2u)+u_r^2/(4u^2)\), direct base
calculation gives
\[
 \operatorname{Ric}_{B,rr}=\overline K-r^{-2}+\frac{u_r}{2ur},
 \qquad
 \operatorname{Ric}_{B,xx}=u\overline K-ur^{-2}+\frac{u_r}{2r},
 \qquad \operatorname{Ric}_{B,rx}=0.
\]
For one fiber, write \(\lambda=\lambda_a\) and
\(\mathsf S=\nabla_B^2\psi_a+d\psi_a\otimes d\psi_a\).
Its entries are
\[
 \mathsf S_{rr}=\lambda_{rr}-\lambda_r/r+\lambda_r^2+r^{-2},\qquad
 \mathsf S_{rx}=\lambda_{rx}-\frac{u_r}{2u}\lambda_x+\lambda_r\lambda_x,
\]
\[
 \mathsf S_{xx}=\lambda_{xx}-(u/r-u_r/2)\lambda_r+u/r^2-u_r/(2r)
                  -\frac{u_x}{2u}\lambda_x+\lambda_x^2.
\]
We use the Einstein equations through order \(r^{2m-2}\).
Subtracting the constant term in the fiber equation gives
\begin{align*}
 0={}&\lambda_{a,rr}-\frac {2N-1}r\lambda_{a,r}
    -\frac{u_r}{2ur}-\frac{\Lambda_r}{r}
    +u^{-1}\lambda_{a,xx}-\frac{u_x}{2u^2}\lambda_{a,x}\\
   &+\frac{u_r}{2u}\lambda_{a,r}
    +\lambda_{a,r}\Lambda_r+u^{-1}\lambda_{a,x}\Lambda_x .
\end{align*}
The \(rr\) and \(xx\) equations at these orders are
\[
 \begin{aligned}
 0={}&\overline K+\frac{u_r}{2ur}-\Lambda_{rr}
              +\frac{\Lambda_r}{r}-\sum_a\lambda_{a,r}^2,\\
 0={}&u\overline K+\frac{Nu_r}{r}-\Lambda_{xx}
       +\left(\frac ur-\frac{u_r}2\right)\Lambda_r
       +\frac{u_x}{2u}\Lambda_x-\sum_a\lambda_{a,x}^2.
 \end{aligned}
\]
The boundary values of \(u^{(1)}\) and \(\Lambda^{(1)}\) vanish,
the latter by \(\sum_a h_a=0\). At order \(r^0\), the linearized
\(xx\) equation and the sum of the linearized fiber equations give
\[
 \begin{aligned}
 (2N-1)[r^2]u^{(1)}+2[r^2]\Lambda^{(1)}&=0,\\
 (2N-1)[r^2]u^{(1)}+2(4N-3)[r^2]\Lambda^{(1)}&=0 .
 \end{aligned}
\]
Thus both \(r^2\) coefficients vanish. For \(2\le q\le m\),
the \(rr\) equation gives
\([r^{2q}](u^{(1)}+2\Lambda^{(1)})=0\). Substituting this into the
\(xx\) equation at order \(r^{2q-2}\) gives
\[
 2q(2q-2N)[r^{2q}]\Lambda^{(1)}
     =\partial_x^2[r^{2q-2}]\Lambda^{(1)} .
\]
Since \(q<N\), induction gives \(u^{(1)}=\Lambda^{(1)}=0\)
through order \(r^{2m}\). The linearized fiber equations now read
\[
 \lambda_{a,rr}^{(1)}-\frac{2N-1}{r}\lambda_{a,r}^{(1)}
           +\lambda_{a,xx}^{(1)}=O(r^{2m}).
\]
The boundary value is \(h_a\cos x/2\), so the radial recurrence
preserves this Fourier mode and gives
\(\lambda_a^{(1)}=h_a A(r)\cos x/2\). Here \(A(0)=1\) and
\[
 2q(2q-2N)[r^{2q}]A=[r^{2q-2}]A,\qquad 1\le q\le m.
\]
Iterating the recurrence and using \(v_a^{(1)}=2\lambda_a^{(1)}\)
gives \eqref{eq:linear-profile}.

The quadratic sources in the \(rr\) and \(xx\) equations therefore
have means
\[
 \sum_a\langle(\lambda_{a,r}^{(1)})^2\rangle=\frac w8(A')^2,
 \qquad
 \sum_a\langle(\lambda_{a,x}^{(1)})^2\rangle=\frac w8A^2.
\]

\subsection{Averaged coefficients of the formal Einstein metric}

\begin{lemma}\label{lem:mean-jets}
Through order \(r^{2m}\), the second means have the form
\[
 \ell_a=-h_a^2/8+wb,\qquad U=w(c-2(2N-1)b),
\]
where \(b,c\) have zero boundary values and coefficients
\[
 c_q=\frac{\mathfrak F_{q-1}}{16q(2N-1)},\qquad
 b_q=\frac{\mathfrak F_{q-1}}{32(2N-1)q(2q-2N)},
 \qquad 1\le q\le m.
\]
Here \(b_q=[r^{2q}]b\) and \(c_q=[r^{2q}]c\). In particular,
\(b_{m}=-c_{m}/(4k)\).
\end{lemma}
\begin{proof}
For \(w=0\) the periodic metric is constant and the second means
vanish, so we may take \(b,c\) with the displayed coefficients.
Suppose \(w>0\). The second-order metric coefficients
enter the equations of Section~\ref{app:coordinate-einstein}
linearly with \(x\)-independent coefficients. Since
\(u^{(1)}=\Lambda^{(1)}=0\), the averaged fiber equations are
\[
 \ell_a''-\frac{2N-1}{r}\ell_a'-\frac{U'}{2r}-\frac{1}{r}\sum_j\ell_j'
   =O(r^{2m}).
\]
Subtract the equations for two fibers. The resulting radial
operator has indicial coefficient \(2q(2q-2N)\ne0\) on \(r^{2q}\)
for \(1\le q\le m\). Thus the terms of positive order in \(\ell_a\)
agree through the determined order. Using the boundary values,
write
\[
 \ell_a=-h_a^2/8+wb,\qquad c=U/w+2(2N-1)b.
\]
Both \(b\) and \(c\) vanish at \(r=0\), and
\(\sum_j\ell_j'=(2N-1)wb'\).

Substitute these expressions into the averaged fiber, \(rr\), and
\(xx\) equations, divide by \(w\), and use
\(\overline K^{(2)}=-u_{rr}^{(2)}/2\). Before division, the quadratic
sources have period averages \(w(A')^2/8\) and \(wA^2/8\).
In the \(rx\) equation, \(\Lambda_{rx}^{(2)}\) averages to zero,
and \(\sum_a\lambda_{a,r}^{(1)}\lambda_{a,x}^{(1)}\) is proportional
to \(\cos x\sin x\). The term \(u_r\Lambda_x/(2u)\) has zero
second coefficient because \(u^{(1)}=\Lambda^{(1)}=0\).
The three averaged equations are
\[
 \begin{aligned}
 b''-\frac{2N-1}{r}b'-\frac{c'}{2r}
    &=O(r^{2m}),\\
 -\frac{c''}{2}+\frac{c'}{2r}-\frac{(A')^2}{8}
    &=O(r^{2m}),\\
 -\frac{c''}{2}+\frac{N}{r}c'
       +(2N-1)\left(b''-\frac{2N-1}{r}b'\right)-\frac{A^2}{8}
    &=O(r^{2m}).
 \end{aligned}
\]
Subtract the second equation and \(2N-1\) times the first from the
third to find \(c'\). Substitution in the first then gives
\begin{equation}\label{eq:mean-system}
 \begin{aligned}
 c'&=\frac{r\mathfrak F}{8(2N-1)}+O(r^{2m+1}),\\
 b''-\frac{2N-1}{r}b'
   &=\frac{\mathfrak F}{16(2N-1)}+O(r^{2m}).
 \end{aligned}
\end{equation}
Comparing coefficients of \(r^{2q-1}\) and \(r^{2q-2}\),
respectively, yields the stated \(c_q\) and \(b_q\). Together with
the boundary values these coefficients determine the mean jets
through order \(r^{2m}\).
\end{proof}

The source in \eqref{eq:mean-system} is \(\mathfrak F=A^2-(A')^2\).
The profile equation and the product rule give
\begin{equation}\label{eq:profile-product}
 (A^2)''=2(A')^2+\frac{2N-1}r(A^2)'+2A^2+O(r^{2m}).
\end{equation}
To compute \(C_j=[r^{2j}]A^2\), use the symmetric-square identity
\[
 \begin{aligned}
 &(A^2)'''-\frac{3(2N-1)}r(A^2)''\\
 &\qquad+\left(\frac{(2N-1)(4N-1)}{r^2}-4\right)(A^2)'
       +\frac{4(2N-1)}r A^2=O(r^{2m-1}).
 \end{aligned}
\]
Comparing the \(r^{2j-3}\) coefficients yields
\begin{equation}\label{eq:profile-convolution}
 \begin{aligned}
 C_j=-\frac{2N-2j+1}{2j(N-j)(2N-j)}C_{j-1},
                         \qquad 1\le j\le m,\\
 C_0=1,\qquad (-4)^m(m!)^2 C_m=2B.
 \end{aligned}
\end{equation}
Multiplying the recurrence gives the identity for \(C_m\) and
\((-1)^jC_j>0\). Combining \eqref{eq:profile-product} with
\(C_{m-1}/C_m=-2mk(N+k)/(2k+1)\) gives
\[
 \begin{aligned}
 \mathfrak F_{m-1}&=2C_{m-1}+2mkC_m
       =-\frac{2mk(2N-1)}{2k+1}C_m,\\
 \frac{4c_m}{C_m}&=-\frac{k}{2(2k+1)}=-a.
 \end{aligned}
\]

\subsection{Proof of the periodic variation formulas}\label{app:variation}

Section~\ref{app:first-profiles} computes \eqref{eq:firstforcing}.
For the second variation, the four terms in \eqref{eq:raw} are obtained
as follows.
\begin{center}
\begin{tabular}{@{}lll@{}}
Term in \eqref{eq:raw} & Coefficient & Calculation \\[2pt]
\(\lvert H^\sharp\eta\rvert^2\) & \(2Y_m=B\) &
  Section~\ref{app:tangential-action} \\
\(\lvert H^\sharp\beta\rvert^2\) & \(2(X_m+\chi/2)=\chi^2-\chi B\) &
  Section~\ref{app:normal-and-norm} \\
\(\lvert H\rvert^2\lvert\eta\rvert^2\) & \(-aB\) &
  \eqref{eq:norm-propagation} \\
\(\lvert H\rvert^2\lvert\beta\rvert^2\) & \(\chi aB\) &
  \eqref{eq:norm-propagation}
\end{tabular}
\end{center}

\subsubsection{The complete second-order expansion}
In the fixed ambient coordinates, \(F_0\) and the boundary contraction
in \eqref{eq:ambient-definition} are independent of \(\varepsilon\).
Write \(\widetilde\Delta_\varepsilon=D_0+\varepsilon D_1+
\varepsilon^2D_2+O(\varepsilon^3)\).
For the constant form \(e\), the ambient codifferential at \(\varepsilon=0\) gives
\(\widetilde\delta_0F_0=-2m(dt/t)\wedge e\), so
\(D_0F_0=0\). Only products ending in a metric variation survive,
with the order of the factors retained:
\[
 [\varepsilon]\widetilde\Delta_\varepsilon^{m}F_0=D_0^{m-1}D_1F_0,
\]
\[
 [\varepsilon^2]\widetilde\Delta_\varepsilon^{m}F_0
   =D_0^{m-1}D_2F_0+
       \sum_{r=1}^{m-1}D_0^{m-1-r}D_1D_0^{r-1}D_1F_0.
\]
In particular, second metric coefficients enter only through \(D_2F_0\);
the sum is empty for \(k=N-1\).
Each differentiated ambient operator preserves closure.
With \(g_0\) as the reference metric, \(W_k(g_0)=\Id\). Since
\(Q_{k,g_0}e=0\), the only additional term in the paired second
coefficient is
\[
 \bigl([\varepsilon]W_k(g_\varepsilon)\bigr)
 \bigl([\varepsilon](Q_{k,g_\varepsilon}e)\bigr).
\]
Flat propagation preserves Fourier modes, so only the period
averages of the second metric coefficients contribute to the mean
of \(D_2F_0\).

Write \(e=dx\wedge\beta+\eta\), with both components transverse.
For an exterior index set \(I\), put \(\sigma=\sum_{i\in I}h_i\)
and \(\kappa=|I|\), so \(\kappa=k\) for components of \(\eta\) and
\(\kappa=k-1\) for components of \(\beta\).
The weighted density for \(dy^I\) is \(t^{\tau_\kappa}\mu\), where
\[
 \tau_\kappa=2N+1-2\kappa,\qquad
 \mu=\sqrt{\det g_r}\prod_{i\in I}v_i^{-1},\qquad v=\mu/u.
\]
\subsubsection{The weighted exterior operators}
Translations and independent reflections of the \(y_a\) coordinates
commute with \(d,\widetilde\delta\). For each \(I\), forms of the type
\(F(t,\rho,x)\wedge dy^I\) have a distinct reflection character.
Both operators preserve these forms. Their pairing is the pairing
induced on \(F\) by the \((t,\rho,x)\) block of \(\widetilde g\),
multiplied by \(\prod_{a\in I}(t^2v_a)^{-1}\).
The \((t,\rho,x)\) block has determinant \(-t^4u\), so the full
density is \(t^{2N+1}\sqrt{\det g_r}\). Multiplying by the \(dy^I\)
pairing factor gives \(t^{\tau_\kappa}\mu\). In the coordinates
\((t,\rho,x)\), the inverse of the ambient metric block is
\[
 (\widetilde g^{ij})_{i,j=t,\rho,x}
 =
 \begin{pmatrix}
 0&-t^{-1}&0\\
 -t^{-1}&2\rho t^{-2}&0\\
 0&0&t^{-2}u^{-1}
 \end{pmatrix}.
\]
Integrating \(d\) against the weighted pairing gives the
codifferential of a \(q\)-form \(F\) in \((t,\rho,x)\) as
\begin{equation}\label{eq:weighted-divergence}
 (\widetilde\delta F)_{i_2\ldots i_q}
 =
 -\frac{\widetilde g_{i_2j_2}\cdots\widetilde g_{i_qj_q}}
        {t^{\tau_\kappa}\mu}
 \partial_i\!\left(t^{\tau_\kappa}\mu F^{i j_2\ldots j_q}\right).
\end{equation}
Here the repeated indices range over \(t,\rho,x\). In the
integration by parts, antisymmetry makes the \(q\) terms from
the exterior derivative identical after relabeling.

For a closed two-form in \((t,\rho,x)\) of dilation degree
\(\omega\), write
\[
 F=t^{\omega-1}(f_1\,dt\wedge d\rho+f_2\,dt\wedge dx)
                         +t^\omega f_3\,d\rho\wedge dx.
\]
Raising indices gives
\[
 F^{t\rho}=-t^{\omega-3}f_1,\qquad
 F^{tx}=-t^{\omega-3}f_3/u,\qquad
 F^{\rho x}=-t^{\omega-4}(f_2-2\rho f_3)/u.
\]
The codifferential has the form
\[
 \widetilde\delta F=t^{\omega-3}D_t\,dt
                   +t^{\omega-2}(D_\rho\,d\rho+D_x\,dx).
\]
For instance, \(\widetilde g_{\rho t}=-t\) gives
\[
 (\widetilde\delta F)_\rho
 =\frac{t}{t^{\tau_\kappa}\mu}
   \bigl\{\partial_\rho(t^{\tau_\kappa}\mu F^{\rho t})
         +\partial_x(t^{\tau_\kappa}\mu F^{xt})\bigr\}
 =t^{\omega-2}\bigl(f_1'+\ell f_1+\mu^{-1}\partial_x(vf_3)\bigr),
\]
where \(\ell=\partial_\rho\log\mu\) and primes denote \(\rho\)-derivatives.
Applying
\eqref{eq:weighted-divergence} to the remaining components gives
\begin{equation}\label{eq:active-adjoint}
 \begin{aligned}
 D_t&=2\rho(f_1'+\ell f_1)-(\tau_\kappa+\omega-3)f_1
                                  +\mu^{-1}\partial_x(vf_2),\\
 D_\rho&=f_1'+\ell f_1+\mu^{-1}\partial_x(vf_3),\\
 D_x&=(\tau_\kappa+\omega-3)f_3
                         +u\mu^{-1}\partial_\rho\{v(f_2-2\rho f_3)\}.
 \end{aligned}
\end{equation}
The next coefficients are
\(f_{1,+}=(\omega-2)D_\rho-\partial_\rho D_t\),
\(f_{2,+}=(\omega-2)D_x-\partial_xD_t\),
\(f_{3,+}=\partial_\rho D_x-\partial_xD_\rho\).
For the normal contribution, \(F=t^{\omega-1}\psi\,dt\wedge d\rho\wedge dx\)
and \(F^{t\rho x}=-t^{\omega-5}\psi/u\). Write
\[
 \widetilde\delta F
   =t^{\omega-3}(f_{1,\delta}\,dt\wedge d\rho
                         +f_{2,\delta}\,dt\wedge dx)
                         +t^{\omega-2}f_{3,\delta}\,d\rho\wedge dx.
\]
Using \eqref{eq:weighted-divergence} with \(q=3\) gives
\[
 f_{1,\delta}=-\mu^{-1}\partial_x(v\psi),\qquad
 f_{3,\delta}=v^{-1}(v\psi)_\rho,
\]
\[
 f_{2,\delta}=2\rho v^{-1}(v\psi)_\rho-(\tau_\kappa+\omega-5)\psi.
\]
Applying \(d\), the top coefficient is
\begin{equation}\label{eq:normal-operator}
 \begin{aligned}
 (\omega-2)f_{3,\delta}-\partial_\rho f_{2,\delta}+\partial_x f_{1,\delta}
  ={}&(\omega-4)v^{-1}(v\psi)_\rho\\
     &-2\rho\partial_\rho\{v^{-1}(v\psi)_\rho\}
       +(\tau_\kappa+\omega-5)\psi_\rho\\
     &-\partial_x\{\mu^{-1}\partial_x(v\psi)\}.
 \end{aligned}
\end{equation}
Write \(f_{j,s}\) and \(\psi_s\) for the coefficients after \(s\)
applications of the ambient Laplacian. At \(t=1,\rho=0\),
the multiplier \(W_k(g_\varepsilon)\) acts
by \(\mu|_{\rho=0}\) on \(dy^I\) and by \(v|_{\rho=0}\) on
\(dx\wedge dy^I\). Thus the paired boundary coefficients are
\(\mu f_{1,m}\) and \(v\psi_m\), evaluated at \(\rho=0\).
The tangential calculation propagates \(\mu f_{1,s}\), which already
includes the metric variation of the pairing. The normal calculation
propagates \(\psi_s\) and adds the pairing variation at the boundary.
Since \(v^{(0)}=1\) and \(\psi_m^{(0)}=0\), the latter gives
\begin{equation}\label{eq:normal-boundary-pairing}
 \left.\langle(v\psi_m)^{(2)}\rangle\right|_{\rho=0}
 =
 \left.\langle\psi_m^{(2)}\rangle\right|_{\rho=0}
 +\left.\langle v^{(1)}\psi_m^{(1)}\rangle\right|_{\rho=0}.
\end{equation}
At step \(s\), both forms have dilation degree \(\omega=2-2s\).

\begin{lemma}\label{lem:finite-jets}
Two diagonal metric families in the periodic model give the same
\(Q_{k,g_\varepsilon}e\) and its first two variations in \(\varepsilon\)
if their \(\rho\)-jets through order \(m\) agree for every \(\varepsilon\).
Thus any smooth extension of the determined metric jet may be used.
\end{lemma}
\begin{proof}
Use the two-form filtration
\(\mathcal I^r=\{f_1,f_3=O(\rho^r),\,f_2=O(\rho^{r+1})\}\).
The \(d\widetilde\delta\) formulas and their metric variations send \(\mathcal I^r\) to \(\mathcal I^{r-1}\).
A metric change \(O(\rho^L)\) changes the operator by a map
\(\mathcal I^r\to\mathcal I^{r+L-1}\). After \(m\) factors
its contribution is in \(\mathcal I^{L-m}\). The normal formula
has the same property on ordinary radial orders. For \(L=m+1\),
the coefficient of \(dt\wedge d\rho\) in the difference of the
final forms is \(O(\rho)\), so it vanishes on restriction to \(\rho=0\).
The same estimates apply after one or two parameter derivatives.
Since \(2m\le 2N-2\), all required metric
coefficients precede the undetermined ambient order.
\end{proof}

\subsubsection{First variation in the periodic model}\label{app:first-profiles}
Put \(\mathsf A(\rho)=A(\sqrt{2\rho})\),
\(\mathsf A_q=\partial_\rho^q\mathsf A\), and
\(K_s=(-2)^{s-1}(k+1)_{s-1}\).
Choose the polynomial \(\mathsf A\) through degree \(m\). Its
equation residual is
\begin{equation}\label{eq:radial-profile}
 \begin{aligned}
 \mathsf R&=2\rho\mathsf A_2-(2N-2)\mathsf A_1-\mathsf A=O(\rho^{m}),\\
 [\rho^q]\mathsf A&=\frac{(-1)^q\Gamma(N-q)}{2^q q!\Gamma(N)},
                       \quad0\le q\le m.
 \end{aligned}
\end{equation}
Write \(\mathsf R_q=\partial_\rho^q\mathsf R\).
Differentiating the radial profile equation gives
\[
 \bigl[-2\rho\partial_\rho^2+2(m-1-2s)\partial_\rho+1\bigr]\mathsf A_q
                      =2(q-k-2s)\mathsf A_{q+1}-\mathsf R_q.
\]
After step \(s\), the remaining \(m-s\) factors each lower the
filtration in Lemma~\ref{lem:finite-jets} by at most one, including
a first metric variation. An error of order \(m-s+1\) thus
vanishes at the boundary and, for \(s<m\), stays above the source
degree \(m-s-1\) at step \(s+1\). These orders refer to
\(\mathcal I^\bullet\) for two-forms and to radial order for the
normal three-form. Multiplication by the pairing density preserves both.

The first variation of \(\log\mu\) and \(\log v\) is
\(-\sigma\mathsf A\cos x\), since \(u^{(1)}=\Lambda^{(1)}=0\) and
\(\lambda_a^{(1)}=h_a\mathsf A\cos x/2\). With \(\psi=1\) and
\(\omega=2\), \eqref{eq:normal-operator} gives the first variation
of the normal coefficient in \(F_1=\widetilde\Delta F_0\):
\[
 \sigma(2\rho\mathsf A_2+2\mathsf A_1-\mathsf A)\cos x
 =\sigma(2N\mathsf A_1+\mathsf R)\cos x.
\]
For the tangential \(D_1F_0\), start with \((f_1,f_2,f_3)=(1,0,0)\) and
\(\omega=2\) in \eqref{eq:active-adjoint}. Then
\(D_t=2\rho\ell-2m\), \(D_\rho=\ell\), and \(D_x=0\).
Taking the first variation of the resulting \(d\widetilde\delta F_0\)
uses \(\ell^{(1)}=-\sigma\mathsf A_1\cos x\) and gives its three
coefficients
\[
 \bigl(2\sigma(\rho\mathsf A_2+\mathsf A_1)\cos x,\,
       -2\rho\sigma\mathsf A_1\sin x,\,
       -\sigma\mathsf A_1\sin x\bigr).
\]

For the induction, use \(V=\mathsf A_s\) with the factor
\(\sigma K_s\) removed. The closed two-form with coefficients
\[
 f_1=2(\rho V'+sV)\cos x,\quad f_2=-2\rho V\sin x,\quad f_3=-V\sin x
\]
agrees with the preceding triple at \(s=1\) after multiplication
by \(\sigma K_1=\sigma\).
At the flat metric, \(\mu=v=u=1\) and
\(\tau_k+\omega-3=2(m-s)\). Substituting the triple into
\eqref{eq:active-adjoint} and using
\(2\rho V''-(2N-2-2s)V'-V=\mathsf R_s\) gives
\[
 \begin{aligned}
 D_t&=\bigl(4(k+s)\rho V'-4s(m-s)V+2\rho\mathsf R_s\bigr)\cos x,\\
 D_\rho&=(2NV'+\mathsf R_s)\cos x,\\
 D_x&=-2(m-s)V\sin x.
 \end{aligned}
\]
Applying \(d\) gives \(-2(k+s)\) times the same coefficient triple
with \(V\) replaced by \(V'\) and \(s\) by \(s+1\), plus the error
\[
 \bigl(-2\{\rho \mathsf R_s'+(s+1)\mathsf R_s\}\cos x,
                      2\rho \mathsf R_s\sin x,\ \mathsf R_s\sin x\bigr).
\]
Since \(\mathsf R_s=O(\rho^{m-s})\), this error lies in
\(\mathcal I^{m-s}\). Induction gives the first variation at step \(s\)
up to \(\mathcal I^{m-s+1}\) for the two-form and
\(O(\rho^{m-s+1})\) for the normal three-form. Also,
\[
 \mathsf A_{s-1}+2N\mathsf A_s
       =2(\rho\mathsf A_{s+1}+s\mathsf A_s)-\mathsf R_{s-1}.
\]
To this accuracy, the first-order coefficients are
\begin{equation}\label{eq:first-profiles}
 \psi_s^{(1)}=2N\sigma K_s\mathsf A_s\cos x,\qquad
 f_{1,s}^{(1)}=\sigma K_s(\mathsf A_{s-1}+2N\mathsf A_s)\cos x,
\end{equation}
\[
 f_{2,s}^{(1)}=-2\rho\sigma K_s\mathsf A_s\sin x,\qquad
 f_{3,s}^{(1)}=-\sigma K_s\mathsf A_s\sin x .
\]
The same estimate discards the initial remainder
\(-\sigma^2\mathsf A \mathsf R/2=O(\rho^{m})\) and the source remainder
\(-N\sigma^2K_s\mathsf A_1\mathsf R_{s-1}=O(\rho^{m-s+1})\) in the
tangential recurrence.

At \(\rho=0\), \(K_{m}\mathsf A_{m}=-1/(2k)\) and
\(K_{m}\mathsf A_{m-1}=1\). Therefore
\[
 \bigl(Q_{k,g_\varepsilon}e\bigr)^{(1)}=-\chi\cos x\,dx\wedge H^\sharp\beta
                     +(1-\chi)\cos x\,H^\sharp\eta .
\]
Since \(Q_{k,g_0}e=0\), the variation of
\(\delta_{g_\varepsilon}Q_{k,g_\varepsilon}e\) is
\(\delta_{g_0}\left.\frac{d}{d\varepsilon}\right|_0
 [W_k(g_\varepsilon)Q_{k,g_\varepsilon}e]\).
Applying the flat codifferential gives
\eqref{eq:firstforcing}.

\subsubsection{The tangential action coefficient}\label{app:tangential-action}
For the tangential contribution, \(\tau_k=2m+1\).
Using \(\mu f_1'=(\mu f_1)'-\ell\mu f_1\) and its derivative,
\eqref{eq:active-adjoint} gives
\begin{align*}
 \mu f_{1,+}={}&-2\rho(\mu f_1)''+(\tau_k+2\omega-7)(\mu f_1)'\\
 &+\ell\{2\rho(\mu f_1)'-(\tau_k+\omega-3)\mu f_1+\partial_x(vf_2)\}\\
 &+(\omega-2)\partial_x(vf_3)-\partial_\rho\partial_x(vf_2).
\end{align*}
The last two terms average to zero. With
\(M=\langle(\log\mu)^{(2)}\rangle\), the initial term
\(\mu f_{1,1}=-2\mu(\rho\ell)'\) gives
\[
 \langle(\mu f_{1,1})^{(2)}\rangle
       =-2(\rho M')'-\frac{\sigma^2}{2}
                           (\mathsf A^2+2N\mathsf A\mathsf A_1)+O(\rho^{m}).
\]
The \(M\)-term is the norm contribution and the \(\sigma^2\)-term
is the action contribution. The first variations are
\[
 \ell^{(1)}=-\sigma\mathsf A_1\cos x,\qquad
 (\mu f_{1,s})^{(1)}=\sigma K_s(\mathsf A_{s-1}+2N\mathsf A_s)\cos x.
\]
Multiplying \(\ell^{(1)}\) by the first variation of its bracket and
averaging over \(x\) gives the source term
\[
 \begin{aligned}
 &\sigma^2K_s\mathsf A_1
     \bigl\{(m-s)(\mathsf A_{s-1}+2N\mathsf A_s)
                    -2N\rho\mathsf A_{s+1}\bigr\}\\
 &\qquad=-\sigma^2K_s\mathsf A_1
          \{(k+s)\mathsf A_{s-1}+2Nk\mathsf A_s\}
          +O(\rho^{m-s+1}).
 \end{aligned}
\]
The reduction uses the differentiated profile equation
\[
 2\rho\mathsf A_{s+1}
   =\mathsf A_{s-1}+(2N-2s)\mathsf A_s+\mathsf R_{s-1}.
\]

The coefficient of \(\rho^{m-s-1}\) in this source depends only on
\[
 T_s=[\rho^{m-s-1}](\mathsf A_1\mathsf A_s),\qquad
 \widehat T_s=[\rho^{m-s-1}](\mathsf A_1\mathsf A_{s-1}),\qquad 1\le s<m.
\]
For the Wronskian
\(\operatorname{Wr}_s=\mathsf A_2\mathsf A_s-\mathsf A_1\mathsf A_{s+1}\), expanding
\(\mathsf R_1\mathsf A_s-\mathsf A_1\mathsf R_s\) from
\eqref{eq:radial-profile} gives
\[
 2\rho\operatorname{Wr}_s'-(2N-3-s)\operatorname{Wr}_s
       =(s-1)(\mathsf A_1\mathsf A_s)'
          +\mathsf R_1\mathsf A_s-\mathsf A_1\mathsf R_s.
\]
For \(s<m-1\), the error has zero coefficient at
\(\rho^{m-s-2}\). The preceding identity gives
\[
 \begin{aligned}
 [\rho^{m-s-2}]\operatorname{Wr}_s&=-\frac{(s-1)(m-s-1)}{2k+s+1}T_s,\\
 T_{s+1}&=\frac{(m-s-1)T_s-[\rho^{m-s-2}]\operatorname{Wr}_s}{2}.
 \end{aligned}
\]
Also \(\widehat T_s=2T_{s+1}-2(N-s)T_s\), by
\eqref{eq:radial-profile}. These identities give
\begin{equation}\label{eq:critical-products}
 \begin{aligned}
 T_{s+1}&=\frac{(k+s)(m-s-1)}{2k+s+1}T_s,\\
 \widehat T_s&=-\frac{2(k+1)(N+k)}{2k+s+1}T_s,
                         \qquad 1\le s<m.
 \end{aligned}
\end{equation}
For \(s=m-1\), the second identity follows directly with \(T_m=0\).
For \(m\ge2\), extracting the \(\rho^{m-1}\) coefficient of
\eqref{eq:profile-product} after \(r^2=2\rho\), and then using
\eqref{eq:profile-convolution}, gives
\[
 \begin{aligned}
 2T_1&=2^{m-1}\bigl(2m(m-N)C_m-C_{m-1}\bigr),\\
 T_1&=\frac{2^{m-1}mk(m-1)}{2k+1}C_m.
 \end{aligned}
\]
For \(m=1\), no \(T_s,\widehat T_s\) occur and all source sums below are empty.

Products of the homogeneous multipliers from step \(s\) to the
boundary give the propagation weights
\[
 p_s=(-2)^{m-s}\frac{(m-1)!(m-s)!}{(s-1)!},
 \qquad p_m=1,\qquad
 p_s=-2s(m-s)p_{s+1}.
\]
Multiplying the recurrence \(z_{s+1}=-2s(m-s)z_s+q_s\) by
\(p_{s+1}\) absorbs the homogeneous multiplier.
Summing \(p_{s+1}z_{s+1}-p_sz_s=p_{s+1}q_s\) gives
\begin{equation}\label{eq:critical-transfer}
 z_m=p_1z_1+\sum_{s=1}^{m-1}p_{s+1}q_s.
\end{equation}
Both action recurrences contain the weighted product
\(p_{s+1}K_sT_s\). With the normalization by \(2s\) below,
the weighted source is a first difference for the tangential
coefficient and for a linear combination of the normal and tangential
coefficients. In both cases the source sum in \eqref{eq:critical-transfer}
telescopes to boundary terms. Define
\[
 \vartheta_1=\frac{Bk}{2m(2k+1)},\qquad
 \vartheta_{s+1}=\frac{(k+s)^2}{s(2k+s+1)}\vartheta_s,
                       \quad 1\le s<m.
\]
By \eqref{eq:profile-convolution}, \(2^mp_1C_m=-B/m^2\).
For \(m\ge2\), the formula for \(T_1\) and
\(p_1=-2(m-1)p_2\) give
\[
 2p_2K_1T_1=-\frac{2^{m-1}p_1C_m\,mk}{2k+1}
             =\frac{Bk}{2m(2k+1)}.
\]
For \(1\le s<m-1\), the ratio of successive terms is
\[
 \begin{aligned}
 \frac{2(s+1)p_{s+2}K_{s+1}T_{s+1}}{2s p_{s+1}K_sT_s}
 &=\frac{s+1}{s}\frac{p_{s+2}}{p_{s+1}}\frac{K_{s+1}}{K_s}
             \frac{T_{s+1}}{T_s}\\
 &=\frac{s+1}{s}\frac{(-2)(k+s)}{-2(s+1)(m-s-1)}
             \frac{(k+s)(m-s-1)}{2k+s+1}\\
 &=\frac{(k+s)^2}{s(2k+s+1)} .
 \end{aligned}
\]
Thus \(2sp_{s+1}K_sT_s=\vartheta_s\) for \(1\le s<m\).
Multiplying the recurrence for \(\vartheta_s\) gives
\[
 \vartheta_m
 =\frac{Bk}{2m(2k+1)}
    \frac{(k+1)_{m-1}^2}{(m-1)!(2k+2)_{m-1}}
 =\frac{1}{2(2k+m)},
\]
where the last equality uses the definition of \(B\).
This formula includes \(m=1\).

For \(1\le s\le m\), let \(Y_s\) be the \(\sigma^2\rho^{m-s}\) coefficient of
\(\langle(\mu f_{1,s})^{(2)}\rangle\), omitting the norm contribution.
Since \([\rho^{m-1}]\mathsf A\mathsf A_1=m2^{m-1}C_m\), its initial value
and recurrence are
\begin{equation}\label{eq:Y-recurrence}
 \begin{aligned}
 Y_1&=-2^{m-2}(C_{m-1}+2NmC_m),\\
 Y_{s+1}&=-2s(m-s)Y_s-K_s\{(k+s)\widehat T_s+2NkT_s\},
                       \qquad 1\le s<m.
 \end{aligned}
\end{equation}
By the recurrence for \(\vartheta_s\), the weighted source is the
first difference of \((s-m)\vartheta_s\),
\[
 -p_{s+1}K_s\{(k+s)\widehat T_s+2NkT_s\}
   =(s+1-m)\vartheta_{s+1}-(s-m)\vartheta_s.
\]
Using \(2^mp_1C_m=-B/m^2\) in the initial coefficient gives
\(p_1Y_1=B/2+(1-m)\vartheta_1\).
Summing these differences from the initial coefficient gives
\begin{equation}\label{eq:Y-endpoint}
 p_sY_s=B/2+(s-m)\vartheta_s,\qquad Y_m=B/2.
\end{equation}
At \(s=m\) the correction vanishes, leaving \(2Y_m=B\) in \eqref{eq:raw}.

\subsubsection{Normal action and norm coefficients}\label{app:normal-and-norm}
The term \(-\partial_x\{\mu^{-1}\partial_x(v\psi)\}\) is an
\(x\)-derivative, so its period mean vanishes. For \(F_0\), use \(\psi=1\) and
\(\omega=2\) in \eqref{eq:normal-operator}. Then
\(v^{-1}(v\psi)_\rho=\partial_\rho\log v\), so the mean normal
coefficient of \(D_2F_0\) is
\[
 \langle\psi_1^{(2)}\rangle
       =-2\partial_\rho\bigl(\rho\partial_\rho
                                \langle(\log v)^{(2)}\rangle\bigr).
\]
This term belongs entirely to the logarithmic density contribution
computed below. The action coefficient \(X_1\) is therefore zero.

For \(v=\mu/u\), the identity \(u^{(1)}=0\) gives
\((\partial_\rho\log v)^{(1)}=-\sigma\mathsf A_1\cos x\).
The radial variation of \eqref{eq:normal-operator} at step \(s\) is
\[
 2\sigma\cos x\,
   \bigl\{\rho\mathsf A_1\partial_\rho+\rho\mathsf A_2
                                  +(s+1)\mathsf A_1\bigr\}.
\]
Its mean on \(\psi_s^{(1)}\) from \eqref{eq:first-profiles} is
\[
 2N\sigma^2K_s\{\rho\partial_\rho+(s+1)\}
                                     (\mathsf A_1\mathsf A_s).
\]

At step \(s\), the radial part of the flat operator on the period
mean is
\[
 -2\rho\partial_\rho^2+2(m-1-2s)\partial_\rho.
\]
On \(\rho^{m-s}\) it gives \(-2s(m-s)\rho^{m-s-1}\).
In the source term, \(\rho\partial_\rho+s+1\) multiplies
\(\rho^{m-s-1}\) by \(m\). Thus the coefficient \(X_s\) of
\(\sigma^2\rho^{m-s}\) in the action term of
\(\langle\psi_s^{(2)}\rangle\), before including the pairing factor \(v\),
satisfies
\begin{equation}\label{eq:X-recurrence}
 X_1=0,\qquad
 X_{s+1}=-2s(m-s)X_s+2NmK_sT_s,\qquad 1\le s<m.
\end{equation}
By \eqref{eq:X-recurrence} and \eqref{eq:Y-recurrence}, the weighted
source for \(X_s+\chi Y_s\) is the first difference of
\(\chi(m+\chi s)\vartheta_s\). Summing gives
\begin{equation}\label{eq:combined-action}
 p_s(X_s+\chi Y_s)
       =\chi(m+\chi s)\vartheta_s,\qquad 1\le s\le m.
\end{equation}
The initial value uses \(p_1Y_1=B/2+(1-m)\vartheta_1\).
At \(s=m\), the identity \((2k+m)\vartheta_m=1/2\) gives
\(X_m+\chi Y_m=\chi(\chi-1)/2\).
At \(\rho=0\), the first variations are
\(v^{(1)}=-\sigma\cos x\) and
\(\psi_m^{(1)}=-\chi\sigma\cos x\), the latter by
\eqref{eq:first-profiles} and \(K_m\mathsf A_m(0)=-1/(2k)\).
The second term of \eqref{eq:normal-boundary-pairing} is therefore
\[
 \left.\langle v^{(1)}\psi_m^{(1)}\rangle\right|_{\rho=0}
 =\chi\sigma^2\langle\cos^2 x\rangle=\frac{\chi\sigma^2}{2}.
\]
The normal paired action coefficient of \(\sigma^2\) is consequently
\[
 X_m+\chi/2=\chi^2/2-\chi Y_m
              =(\chi^2-\chi B)/2.
\]

The norm terms come from the second logarithmic means of the
densities. For \(|I|=\kappa\), Lemma~\ref{lem:mean-jets} gives
\[
 \langle(\log\mu)^{(2)}\rangle
    =U/2+(2N-1-2\kappa)wb-w/8+\tfrac14\sum_{a\in I}h_a^2.
\]
For \(v=\mu/u\), the same formula has \(-U/2\) in place of
\(U/2\). Substitute \(U=w(c-2(2N-1)b)\), with \(\kappa=k\)
for the tangential contribution and \(\kappa=k-1\) for the normal
contribution. Modulo terms constant in \(r\), this gives
\[
 \langle(\log\mu)^{(2)}\rangle=w(c/2-2kb),\qquad
 \langle(\log v)^{(2)}\rangle=w(-c/2+2(2N-k)b).
\]
For either logarithmic mean \(M\), the initial norm term is
\(-2(\rho M')'\). Subsequent norm coefficients satisfy the homogeneous
recurrence with multiplier \(-2s(m-s)\), so \eqref{eq:critical-transfer}
gives
\begin{equation}\label{eq:norm-propagation}
 p_1[\rho^{m-1}]\bigl(-2(\rho M')'\bigr)
 =-2m^2p_1[\rho^m]M
 =(-2)^m(m!)^2[\rho^m]M.
\end{equation}
Since \(r^2=2\rho\) and \(b_m=-c_m/(4k)\), the coefficients
\([\rho^m]M\) are \(2^mwc_m\) for the tangential contribution and
\(-2^m\chi wc_m\) for the normal contribution.
By \eqref{eq:profile-convolution} and \(4c_m/C_m=-a\),
\eqref{eq:norm-propagation} yields the norm coefficients
\(-aBw/2\) and \(\chi aBw/2\) in the paired \(\varepsilon^2\) term.

\begin{proof}[Proof of Proposition~\ref{prop:raw-hessian}]
The calculations above give \eqref{eq:firstforcing} and the four
coefficients in \eqref{eq:raw} on exterior basis forms.
The distinct fiber-reflection characters exclude cross terms,
so summing over the basis gives \eqref{eq:raw} for general \(e\).
\end{proof}

\subsection{Checks of the variation formulas}\label{sec:variation-checks}

\subsubsection{The case \(m=1\)}
When $m=1$, the formula of
\cite[Proposition 8.4]{AG}, normalized to unit
leading coefficient, takes the form
\[
 Q_{k,g}e=\Delta_ge-\frac{2}{k}(g^{-1}\operatorname{Ric}_g)^\sharp e
              +\frac{2}{2k+1}\operatorname{Scal}_g e .
\]
Here the action of the Ricci endomorphism on forms is the same
slotwise action as for \(H^\sharp\). For diagonal \(H\) it suffices
to take \(e=dy^I\) or \(e=dx\wedge dy^I\). Set
\(\sigma=\sum_{i\in I}h_i\) and \(w=|H|^2\).
The norm density \(|e|^2_{g_\varepsilon}dV_{g_\varepsilon}/dV_{g_0}\)
has first Taylor coefficient \(-\sigma\cos x\) in both cases.
To take paired second coefficients, pair each term with \(e\) in
\(g_\varepsilon\), multiply by \(dV_{g_\varepsilon}/dV_{g_0}\),
and extract \([\varepsilon^2]\).
With \(p_a=\tfrac12\log(1+\varepsilon h_a\cos x)\), the Ricci eigenvalues
in the \(x\) and \(y_a\) directions are
\(-\sum_a(\partial_x^2p_a+(\partial_xp_a)^2)\) and
\(-\partial_x^2p_a-(\partial_xp_a)\sum_b\partial_xp_b\), respectively.
Thus the mean second coefficient of the paired
Ricci term is \(-\sigma^2/4\) for \(dy^I\) and
\(-\sigma^2/4-w/8\) for \(dx\wedge dy^I\). The mean second coefficient
of the paired scalar curvature term before its prefactor is
\(-w/8\) in both cases. Since \(de=0\), the Laplacian pairing is
\(\|\delta_{g_\varepsilon}e\|^2_{g_\varepsilon}\). It contributes
zero for \(dy^I\) and has mean second coefficient \(\sigma^2/2\)
for \(dx\wedge dy^I\), for which
\((\delta_ge)^{(1)}=-\sigma\sin x\,dy^I\).
The first coefficient of \(Q_{k,g}e\) is
\(-(\sigma/k)\cos x\,dy^I\) or
\(-(1+1/k)\sigma\cos x\,dx\wedge dy^I\), respectively.
Applying the flat codifferential gives \eqref{eq:firstforcing}.
Doubling the paired second coefficients gives
\[
 \begin{aligned}
 R(dy^I;dx,H)&=\frac{\sigma^2}{k}-\frac{w}{2(2k+1)},\\
 R(dx\wedge dy^I;dx,H)&=
   \left(1+\frac1k\right)\sigma^2+\frac{(k+1)w}{2k(2k+1)},
 \end{aligned}
\]
as in \eqref{eq:raw} with \(B=1/k\).

\subsubsection{The six-dimensional relative problem}\label{app:six-dimensional-check}
For \((N,k)=(3,1)\), the relative problem of
\cite[Proposition 3.11 and Section 4.2.3]{AG} gives a separate check.
The closed boundary one-form \(e\) is the leading normal coefficient of a
formal harmonic two-form with zero tangential trace. With \(dr/r\) placed
first in the leading and logarithmic terms, \(Q_{1,g}e\) is \(-16\)
times the \(r^4\log r\) coefficient of the normal component
\cite[Theorem 6.2(ii)]{AG}.
Write \(\mu_I=\sqrt{\det g_r}\prod_{i\in I}v_i^{-1}\), and for a radial
density \(q\) set \(\ell_j(q)=[r^j]\log q\) and \(q_0=q|_{r=0}\).

For \(e=dy^a\), write the closed extension as
\(dr/r\wedge dy^a+d(\varphi\,dy^a)\), where
\(\varphi=r^2\varphi_2+r^4(\varphi_{\log}\log r+\varphi_4)+\cdots\).
Its codifferential must be \(4\,dy^a\):
closedness of the codifferential makes its scalar coefficient
constant, and the leading term fixes that constant.
The equation for \(\varphi\), with \(\mu=\mu_{\{a\}}\), is
\[
 \bigl((r\partial_r)^2-4r\partial_r\bigr)\varphi
   +(r\partial_r\log\mu)(1+r\partial_r\varphi)
   +r^2\mu^{-1}\partial_x\bigl((\mu/u)\partial_x\varphi\bigr)=0.
\]
At order \(r^2\) this gives \(\varphi_2=\ell_2(\mu)/2\). At order \(r^4\)
the logarithmic term contributes \(4\varphi_{\log}\). Multiplication by
\(\mu_0\) and integration over the period remove the \(x\)-divergence
since \(u|_{r=0}=1\).
The paired mean of the normal logarithmic coefficient is therefore
\[
 \langle\mu_0(4\varphi_{\log})\rangle
    =-4\langle\mu_0\{\ell_4(\mu)+\ell_2(\mu)\varphi_2\}\rangle.
\]
By \eqref{eq:linear-profile} and Lemma~\ref{lem:mean-jets},
\[
 \begin{gathered}
 \mu_0^{(1)}=-h_a\cos x,\qquad
 \ell_2(\mu)^{(1)}=\frac{h_a}{8}\cos x, \\
 \ell_4(\mu)^{(1)}=-\frac{h_a}{64}\cos x,\qquad
 \langle\ell_4(\mu)^{(2)}\rangle=-\frac{w}{512}.
 \end{gathered}
\]
Thus the mean second coefficient is \(w/128-3h_a^2/64\).
The \(r^4\) equation before averaging gives
\((4\varphi_{\log})^{(1)}=h_a\cos x/8\).

For \(e=dx\), the closed extension is \(p\,dr/r\wedge dx\), with
\(p=1+r^2p_2+r^4(p_{\log}\log r+p_4)+\cdots\).
Writing \(v=\mu_{\varnothing}/u\), the condition that its codifferential
is closed becomes
\[
 \partial_r\{4p-r\partial_r p-r(\partial_r\log v)p\}
    =\partial_x\{r\mu_{\varnothing}^{-1}\partial_x(vp)\}.
\]
At order \(r\), this gives
\(4(p_2-\ell_2(v))=\partial_x(\mu_{\varnothing,0}^{-1}\partial_x v_0)\).
The first coefficients of \(p_2,\ell_2(v),\ell_4(v)\), and \(v_0\)
all vanish by trace-freeness and \(u^{(1)}=0\).
At order \(r^3\), taking the mean gives
\(\langle p_{\log}\rangle
=-2\langle\ell_2(v)p_2\rangle-4\langle\ell_4(v)\rangle\).
Here \(\langle\ell_4(v)^{(2)}\rangle=3w/512\), by
Lemma~\ref{lem:mean-jets}. Hence the paired mean second coefficient
of \(p_{\log}\) is \(-3w/128\), and \(p_{\log}^{(1)}=0\).
Multiplication by \(-16\) gives first variations of \(Q_{1,g}e\)
that are coclosed at the flat metric for both types of input, in agreement
with \eqref{eq:firstforcing} since \(H^\sharp\beta=0\) for \(k=1\).
Multiplication by \(-32\) gives
\[
 R(dy^a;dx,H)=\frac32h_a^2-\frac{w}{4},\qquad
 R(dx;dx,H)=\frac{3w}{4}.
\]
These agree with \eqref{eq:raw} for \(B=3/2\), \(\chi=3\), and \(a=1/6\).

\end{document}